\documentclass[11pt]{article}
\usepackage[centertags]{amsmath}
\usepackage{amsfonts}
\usepackage{amssymb}
\usepackage{amsthm}
\usepackage{newlfont}
\usepackage{stmaryrd}
\usepackage{color}

\theoremstyle{plain}
\newtheorem{thm}{Theorem}[section]
\newtheorem{cor}[thm]{Corollary}
\newtheorem{lem}[thm]{Lemma}
\newtheorem{prop}[thm]{Proposition}
\newtheorem{exam}[thm]{Example}
\newtheorem{rem}[thm]{Remark}
\newtheorem{defi}[thm]{Definition}
\newtheorem{ques}[thm]{Question}
\newtheorem{conj}[thm]{Conjecture}

\def\sqr#1#2{{\vcenter{\vbox{\hrule height.#2pt
              \hbox{\vrule width.#2pt height#1pt \kern#1pt \vrule
width.#2pt}
              \hrule height.#2pt}}}}
\def\3n{\negthinspace \negthinspace \negthinspace }
\def\2n{\negthinspace \negthinspace }
\def\1n{\negthinspace }

\def\no{\noindent}

\def\ms{\medskip}
\def\bs{\bigskip}

\def\be{\begin{equation}}
\def\bel{\begin{equation}\label}
\def\ee{\end{equation}}
\def\bea{\begin{eqnarray}}
\def\eea{\end{eqnarray}}
\def\bt{\begin{theorem}}
\def\et{\end{theorem}}
\def\bc{\begin{corollary}}
\def\ec{\end{corollary}}
\def\bl{\begin{lemma}}
\def\el{\end{lemma}}
\def\bp{\begin{proposition}}
\def\ep{\end{proposition}}
\def\br{\begin{remark}}
\def\er{\end{remark}}
\def\ba{\begin{array}}
\def\ea{\end{array}}
\def\bd{\begin{definition}}
\def\ed{\end{definition}}

\newcommand{\poly}{\mathbb{C}[z_1,\ldots,z_d]}

\newcommand{\Id}{\mathrm{Id}}

\newcommand{\card}{\mathrm{card~}}

\newcommand{\range}{\mathrm{ran~}}
\newcommand{\coker}{\mathrm{coker~}}

\newcommand{\Tr}{\operatorname{Tr}}
\newcommand{\Ext}{\mathrm{Ext}}
\newcommand{\ord}{\mathrm{ord}}
\newcommand{\ind}{\mathrm{ind}}

\makeatletter
   
   \@addtoreset{equation}{section}
\makeatother

\begin{document}

\title{\bf Essential normality of quotient modules vs. Hilbert-Schmidtness of submodules in $H^2(\mathbb D^2)$}
\author{Penghui Wang\thanks{Partially supported by NSFC(No. 12271298,12671156) and Taishan Scholars Project (TSTP 20250702), E-mail: phwang@sdu.edu.cn. \ms}
\quad Chong Zhao \thanks{Partially supported by NSFC(No. 12071253,12571137), and the Young Scholars Program of Shandong University, {E-mail:} { chong.zhao@sdu.edu.cn}. \ms}\quad Zeyou Zhu\thanks {Partially supported by NSFC(No. 11871308,12271298), E-mail: 13155486329@163.com \ms}
 \\ \\
 School of Mathematics, Shandong University\\
Jinan, Shandong 250100, The People's Republic of China\\}

\maketitle
\begin{abstract}

In the present paper, the essential normality of non-algebraic quotient modules over the polydisc is considered. We prove that any quotient module in an integer-weighted Bergman module over the bidisk, associated to a finitely generated submodule containing a certain distinguished polynomial, is essentially normal. As an application, we prove that each finitely generated submodule containing a nonzero polynomial is Hilbert-Schmidt, which partially solves the conjecture of Yang \cite{Ya3}.
\end{abstract}

\bs

\no{\bf 2020 Mathematics Subject Classification}. 47A13, 46H25

\bs

\no{\bf Key Words}. Distinguished variety, essential normality, non-algebraic quotient module, Hilbert-Schmidt submodule, weighted Bergman space over the polydisc

\section{Introduction}

~~~~In the present paper, we consider the essential normality of non-algebraic quotient modules in the weighted Bergman modules over the unit polydisc $\mathbb{D}^d$ for $d>1$. For a real number $\alpha>-1$, the weighted Bergman space over $\mathbb{D}^d$ is defined by
$$
{\cal A}_\alpha(\mathbb D^d)=\left\{f\in\text{Hol}(\mathbb{D}^d):\|f\|_\alpha^2=c_\alpha\int_{\mathbb D^d}|f(z)|^2\prod\limits_{i=1}^d(1-|z_i|^2)^\alpha \text{dV}_d(z)<\infty\right\},
$$
where $\text{dV}_d$ is the usual Lebesgue measure on $\mathbb D^d$ and $c_\alpha$ is the constant that normalizes the measure $\prod\limits_{i=1}^d(1-|z_i|^2)^\alpha\text{dV}_d$. It can be verified that ${\cal A}_\alpha(\mathbb D^d)$ is a Hilbert space with reproducing kernel
\begin{equation}
K_{\lambda,\alpha}(z)=\prod_{i=1}^d(1-\bar{\lambda}_iz_i)^{-\alpha-2},\quad \lambda,z\in\mathbb D^d.
\end{equation}
The weighted Bergman space ${\cal A}_\alpha(\mathbb D^d)$ can be viewed as a Hilbert module \cite{CG,DP} over the polynomial ring $\poly$, endowed with a natural module action defined by the multiplication by polynomials. A closed subspace $\cal M$ of ${\cal A}_\alpha(\mathbb D^d)$ that is invariant under the module action is called a submodule, and $\cal N={\cal M}^\perp$ is called the associated quotient module. The module actions on submodules and quotient modules are given by
 $$
 p\cdot h=R_{p,\mathcal{M}}h,~p\in\poly,~h\in \cal M
 $$
 and
 $$
 p\cdot g=S_{p,\mathcal{N}}g,~p\in\poly,~g\in\cal N
 $$
respectively,  where $R_{p,\mathcal{M}}=M_p|_{\cal M}$ and $S_{p,\mathcal{N}}=P_{\cal N}M_p|_{\cal N}$, and $P_{\cal N}$ is the orthogonal projection onto $\cal N$. For brevity, we write $$R_{p}=R_{p,\mathcal{M}} \text{ and }S_{p}=S_{p,\mathcal{N}}$$
if there is no confusion on the modules being referred to. A quotient module $\cal N$ of ${\cal A}_\alpha(\mathbb D^d)$ is said to be essentially normal if all the commutators $[S_{p}^*,S_{p}]$ for $p\in\poly$ are compact. Let $\cal K(\cal N)$ be the ideal of compact operators on $\cal N$, written as $\cal K$ for short. When $\alpha=-1$, $K_{\lambda,\alpha}$ is just the reproducing kernel of the classical Hardy space $H^2(\mathbb D^d)$, hence we also use ${\cal A}_{-1}(\mathbb D^d)$ to refer to $H^2(\mathbb D^d)$. For brevity, the closure of an ideal $I\subset\poly$ in $\mathcal{A}_\alpha(\mathbb D^d)$ is denoted $[I]_\alpha$. We also write $[f_1,\ldots,f_n]_\alpha$ for the submodule of $\mathcal{A}_\alpha(\mathbb D^d)$ generated by functions $f_1,\ldots,f_n$.

In \cite{Arv2}, Arveson conjectured that all homogeneous submodules of the symmetric Fock module over the complex unit ball are essentially normal. Much work has been done along this line, such as \cite{Arv2,Dou1,Dou2,DTY,DW,DWY,Guo1,GWk2,GWk1,KS,Sha,WX3} and the references therein, and \cite{GWy} is a good survey on this problem. A polydisc version of Arveson's conjecture is to characterize the essential normality of homogeneous quotient modules in $H^2(\mathbb D^d)$, which was solved by the first two authors in \cite{WZ2}. Other efforts on the polydisc version of the Arveson's conjecture can be found in \cite{Cla,DGS,GWp1,GWp2,Wa,WY,WZ1,WZ3} and references therein.

It is known that the essential normality of quotient modules of analytic Hilbert modules over the polydisc is closely related to distinguished varieties, which were introduced by Agler and $\operatorname {M^cCarthy}$ \cite{AM1} to deal with Ando's inequality and Pick interpolation problem, etc. An ideal $I$ in $\poly$ is said to be distinguished if each component of the variety $Z(I)$ intersects $\mathbb{D}^d$ and satisfies
$$\partial\big(Z(I)\cap\mathbb D^d\big)\subset \mathbb T^d.$$
In this case, the submodule $[I]_\alpha\subset{\cal A}_\alpha(\mathbb D^d)$ generated by $I$ and the associated quotient module $[I]_\alpha^\bot$ are said to be distinguished. As usual, we write $(z,w)$ for a $2$-dimensional variable. A polynomial $p\in\mathbb{C}[z,w]$ is said to be distinguished [1] if the principal ideal generated by $p$ is distinguished. For the essential normality of inhomogeneous distinguished quotient weighted Bergman modules, recently in \cite{GWZ}, by introducing the Grassmannian structure for the distinguished quotient modules, Guo and the first two authors proved the following result.
  \begin{thm}(\cite[Theorem 1.3]{GWZ})\label{them:GWZ}
  Let $V$ be an algebraic subvariety of $\mathbb D^2$, and $I_V$ be the radical ideal of $\mathbb C[z,w]$ determined by $V$. Then $V$ is distinguished if and only if $[I_V^2]_{\alpha}^\perp$ is essentially normal.

Furthermore, let $I\subset \mathbb C[z,w]$ be a distinguished ideal, then for every integer $$\alpha\geq\max\{m_{Z(I)\cap\mathbb D^2},n_{Z(I)\cap\mathbb D^2}\}-2,$$ the quotient module $[I]_\alpha^\perp$ is $(1,\infty)$-essentially normal.

  \end{thm}
For a given distinguished variety $V\subset\mathbb{D}^2$, set
$$m_V=\max_{\lambda\in\mathbb{D}}\card\{w\in\mathbb{D}:(\lambda,w)\in V\},~n_V=\max_{\lambda\in\mathbb{D}}\card\{z\in\mathbb{D}:(z,\lambda)\in V\}.$$
In the language of Agler and M\raisebox{.5ex}{c}Carthy \cite{AM1}, $(m_V,n_V)$ is called the rank of $V$. In what follows, for a distinguished variety $V$ in $\mathbb D^2$, $p_V$ denotes a distinguished polynomial with $Z(p_V)\cap\mathbb D^2=V.$ The following result characterizes the essential normality of non-algebraic quotient modules associated to distinguished varieties.
\begin{thm}\label{Theorem:Main}
Let $V$ be a distinguished variety and $\alpha\geq\max\{m_{V},n_{V}\}-2$ be an integer. Then for every finitely generated submodule ${\cal M}_{V,\alpha}$ of ${\cal A}_\alpha(\mathbb D^2)$ containing a distinguished polynomial $p_V$, the quotient module ${\cal M}_{V,\alpha}^\perp$ is $(1,\infty)$-essentially normal.
\end{thm}
It is worth pointing out that Theorem \ref{Theorem:Main} is the first result on the essential normality of non-algebraic quotient modules over the polydisc, where a quotient module is said to be non-algebraic if the associated submodule cannot be generated by polynomials. In this theorem, the choice of weighted Bergman space depends on the algebraic rank of the variety. Therefore, it is worth finding a fixed analytic Hilbert module over $\mathbb D^2$, such as the Hardy space $H^2(\mathbb D^2)$, on which all the distinguished quotient modules are essentially normal. There are two important results along this line. Firstly, it was proved in \cite{WZ1,WZ3} that all the distinguished (quasi)-homogeneous quotient modules of $H^2(\mathbb D^d)$ are essentially normal. Secondly, Clark \cite{Cla} proved that for nonconstant finite Blaschke products $B_{i}(z_i)$, $i=1,\ldots,d$, the quotient Hardy module $[B_1(z_1)-B_2(z_2),\cdots,B_{d-1}(z_{d-1})-B_d(z_d)]^\perp$ is a variant of weighted Bergman space over the distinguished variety $\mathbb{D}^d\cap\bigcap\limits_{i=1}^{d-1} Z\big(B_{i}(z_i)-B_{i+1}(z_{i+1})\big)$, hence is essentially normal. Let $\eta_1, \eta_2$ be two nonconstant inner functions, by \cite{Cla} and \cite{Wa}, $[\eta_1(z)-\eta_2(w)]^\perp$ is essentially normal if and only if both $\eta_1$ and $\eta_2$ are finite Blaschke products. For higher dimensional cases, we have the following result, which is a special case of Corollary \ref{cor:ess}.
\begin{prop}

Let $\eta_1,\ldots,\eta_d$ be nonconstant inner functions and $n_1,\ldots,n_{d-1}$ be positive integers. Then the quotient Hardy module $\big[\big(\eta_1(z_1)-\eta_2(z_2)\big)^{n_1},\cdots,\big(\eta_{d-1}(z_{d-1})-\eta_d(z_d)\big)^{n_{d-1}}\big]^\perp$ is essentially normal if and only if all the $\eta_i$'s are finite Blaschke products.
\end{prop}

Using the notations in Clark \cite{Cla} and Ferguson and Rochberg \cite{FR1},  the quotient Hardy module $\big[\big(\eta_1(z_1)-\eta_2(z_2)\big)^{n_1},\cdots,\big(\eta_{d-1}(z_{d-1})-\eta_d(z_d)\big)^{n_d}\big]^\perp$ can be seen as the higher order restriction of $H^2(\mathbb D^d)$ on the distinguished variety $\mathbb{D}^d\cap\bigcap\limits_{i=1}^{d-1} Z\big(B_{i}(z_i)-B_{i+1}(z_{i+1})\big)$, which is essentially normal.
 However, even the essential normality of $\big[\big(B_1(z)-B_2(w)\big)^N\big]^\bot~(N\geq3)$ has not been established before. Inspired by the above observations, we introduce the following notations. For a tuple of non-constant finite Blaschke products $B=\big(B_1(z_1),B_2(z_2),\cdots,B_d(z_d)\big)$ and a polynomial $p$, write $p_B(z)=p\big(B_1(z_1),\cdots,B_d(z_d)\big).$

The following result gives  the essential normality of non-algebraic quotient modules in $H^2(\mathbb D^2)$.
\begin{thm}\label{Theorem:Main1}
    Let $p$ be a distinguished quasi-homogeneous polynomial in $\mathbb C[z,w]$, and $p_B$ be defined as above. Then for any finitely generated Hardy submodule $\mathcal{M}\subset H^2(\mathbb{D}^2)$ containing $p_B$, the quotient $\mathcal{M}^\perp$ is essentially normal.
\end{thm}

\begin{rem}
In Theorem \ref{Theorem:Main1}, the assumption that $p$ is distinguished  cannot be omitted in general. In fact, assuming that the order of $B_1(z_1)$ is at least $2$, and $p$ has at least one root on $\mathbb{D}\times\mathbb{T}$, then an application of Berezin transform shows that $[p_{B}]^\perp$ is not essentially normal.
\end{rem}

 More surprisingly, in the case that a submodule $[I]$ contains a $p_B$ as in Theorem \ref{Theorem:Main1}, we find that the essential normality of $[I]^\perp$ is equivalent to the Hilbert-Schmidtness of $[I]$. To state this more clearly, we recall some concepts. For a submodule $\mathcal{M}$ in $H^2(\mathbb D^2)$, the core operator
$$
C_\mathcal{M}=\Id_{\mathcal{M}}-R_{z}R_{z}^*-R_{w}R_{w}^*+R_{zw}R_{zw}^*,
$$
introduced by Guo and Yang \cite{GY} plays an important role in the study, where $\Id_{\mathcal{M}}$ is the identity operator on $\mathcal{M}$. A submodule ${\mathcal{M}}$ is said to be Hilbert-Schmidt if $C_\mathcal{M}$ is a Hilbert-Schmidt operator. Yang \cite[page 403]{Ya3} posed the following conjecture.
\begin{conj}\label{Conjecture:HS}
All finitely generated submodules of $H^2(\mathbb D^2)$ are Hilbert-Schmidt.
\end{conj}
Some partial positive results on this conjecture have been obtained. In \cite{Ya3,Ya2}, it was proved that all algebraic submodules of $H^2(\mathbb{D}^2)$ are Hilbert-Schmidt. Luo, Izuchi and Yang \cite{LIY} proved that a submodule of $H^2(\mathbb D^2)$ containing $z-B(w)$ is Hilbert-Schmidt if and only if it is finitely generated, where $B$ is a nonconstant finite Blaschke product. This result was then generalized by Zu, Yang and Lu \cite{ZYL} to the submodules containing $B_1(z)-B_2(w)$ for finite Blaschke products $B_1$ and $B_2$.
As an application of the idea in the proof of Theorem \ref{Theorem:Main1}, we have the following result, which addresses Conjecture \ref{Conjecture:HS}.
\begin{thm}\label{thm:main3}
Let $\mathcal{M}$ be a finitely generated submodule of $H^2(\mathbb{D}^2)$ containing a nonzero polynomial, then $\mathcal{M}$ is Hilbert-Schmidt.
\end{thm}

At last, after proving the Hilbert-Schmidtness of submodules, one will be interested in the calculation of Hilbert-Schmidt norms of the core operator $C_{\cal M}$ for some Hilbert-Schmidt submodules $\cal M$. The toy model is the submodule $\mathcal M_{N;B_1,B_2}=[(B_1(z)-B_2(w))^N]$ for finite Blaschke products $B_i$, which is the simplest model considered in the present paper. It is surprising to find that there is no difference between the trace-type formulas for $\mathcal{M}_{N;B_1,B_2}$ and the formula for $\mathcal{M}_{N;\theta,\psi}$, where $\theta$ and $\psi$ are nonconstant inner functions. It will be seen that the Schatten-$p$ norm depends only on $N$ and the pseudohyperbolic distance between $\theta(0)$ and $\psi(0)$. Since the trace formula is a little technical, we put it in Appendix.

This paper is organized as follows. In Section 2, we obtain the essential normality of the so-called Clark-type quotient modules of $H^2(\mathbb{D}^d)$ which is the foundation for proving Theorem \ref{Theorem:Main1}. In Section 3, the K-homology elements of Clark-type quotient modules are considered, and the non-triviality of the K-homology elements induced by the associated $C^*$-extensions is obtained. In Section 4, for submodules $\mathcal{M}$ containing a Clark-type ideal, the equivalence between the Hilbert-Schmidtness of $\mathcal{M}$ and the essential normality of $\mathcal{M}^\perp$ is established; moreover, Theorems \ref{Theorem:Main} and \ref{Theorem:Main1} will be proved. In Section 5, finitely generated non-algebraic submodules of $H^2(\mathbb{D}^2)$ containing a nonzero polynomial are shown to be Hilbert-Schmidt.  In Section 6, several examples on non-algebraic submodules and the associated quotient modules are investigated. Finally, we will derive some interesting trace formulas in Appendix.
\section{$(1,\infty)$-essentially normal quotient modules}

~~~~~Recall that Clark's quotient module $[B_1(z_1)-B_2(z_2),\cdots, B_{d-1}(z_{d-1})-B_d(z_d)]^\perp$ was proven essentially normal in \cite{Cla}. Inspired by this fact, we generalize Clark's quotient module to the concept of $\emph{Clark-type quotient module}.$  To simplify the notations, for $q\in\poly$ and a $d$-tuple of inner functions $\eta(z)=\big(\eta_1(z_1),\ldots,\eta_d(z_d)\big),$ we write
$$q_\eta(z)=q\big(\eta_1(z_1),\ldots,\eta_d(z_d)\big),$$
and $J_\eta=\{q_\eta:q\in J\}$ for a subset $J\subset\poly$. If $[J]^\perp$ is of finite dimension, then obviously, it is essentially normal. \emph{In what follows, we always assume that $[J]^\perp$ is of infinite dimension.  }

To continue, we need the terminology of quasi-homogeneity. By the term weight we mean a given tuple $\mathbf{K}=(K_1,\ldots,K_d)$ of positive integers. For a monomial $z^\alpha$, the integer $\alpha\cdot\mathbf{K}=\sum_{i=1}^d\alpha_iK_i$ is called its $\mathbf{K}$-degree. A linear combination of monomials with the same $\mathbf{K}$-degree is said to be $\mathbf{K}$-quasi-homogeneous. Ideals generated by $\mathbf{K}$-quasi-homogeneous polynomials are said to be $\mathbf{K}$-quasi-homogeneous, and so are the associated submodules and quotient modules. In \cite{WZ3}, the essential normality of $\mathbf{K}$-quasi-homogeneous quotient modules of $H^2(\mathbb D^d)$ was completely characterized. In particular, distinguished $\mathbf{K}$-quasi-homogeneous quotient Hardy modules are essentially normal.

\begin{defi}\label{def:Clarktype}
Let $I\subset\poly$ be an ideal. If there exists a $d$-tuple of nonconstant finite Blaschke products $B(z)=\big(B_1(z_1),\ldots,B_d(z_d)\big)$ and a distinguished quasi-homogeneous ideal $J$ of $\poly,$ such that
\begin{eqnarray*}
    I&=&\big(J_B\poly\big)\cap\poly,
\end{eqnarray*}
then we call $I$ a Clark-type ideal, $Z(I)\cap\mathbb{D}^d$ a Clark-type variety, and $[I]^\perp$ a Clark-type quotient module.
\end{defi}
From the fact that $|z_i|\to1$ if and only if $|B_i(z_i)|\to1~(1\leq i\leq d)$, we see that Clark-type varieties are distinguished. Since the denominators of finite Blaschke products are invertible outer functions, the submodule $[I]$ and the submodule generated by $J_B$ coincide. Note that an ideal with Clark-type variety is not necessarily Clark-type.


To precisely characterize the essential normality of quotient modules over $\mathbb{D}^d$, we introduce the concept of Macaev ideals. Let $\lambda_1,\lambda_2,\ldots$ be the list of positive singular values of a compact operator $A\in B(H)$ in the descending order, repeated according to multiplicity. If the norm $\|A\|_{1,\infty}=\sup_{n\geq1}\big|\frac{1}{\ln(n+1)}{\sum\limits_{k=1}^n\lambda_k}\big|$ is finite, then $A$ is said to be in the Macaev ideal $\mathfrak{L}^{(1,\infty)}$. The class $\mathfrak{L}^{(2,\infty)}$ is defined as the collection of operators $A\in B(H)$ with $A^*A\in\mathfrak{L}^{(1,\infty)}$, with norm $\|A\|_{2,\infty}=\sqrt{\|A^*A\|_{1,\infty}}$. For $A\in\mathfrak{L}^{(2,\infty)}$ and $B\in B(H)$, it holds that $\|A^*\|_{2,\infty}=\|A\|_{2,\infty}$ and $\|BA\|_{2,\infty}\leq\|B\|\cdot\|A\|_{2,\infty}$. The reader is referred to \cite{Conn} for more details of Macaev ideals. As in \cite{WZ1}, a quotient module $N$ in $H^2(\mathbb{D}^d)$ is said to be $(1,\infty)$-essentially normal if $[S_p^*,S_q]\in\mathfrak{L}^{(1,\infty)}$ for all polynomials $p$ and $q$.

For $A\in B(H)$, if it happens that $\Id_H-A^*A\in\mathfrak{L}^{(1,\infty)}$, then $A$ is said to be $(1,\infty)$-essentially isometric. The adjoint operators of $(1,\infty)$-essential isometries are called $(1,\infty)$-essential co-isometries. If an operator is simultaneously $(1,\infty)$-essentially isometric and $(1,\infty)$-essentially co-isometric, then it is said to be $(1,\infty)$-essentially unitary.

To prove the essential normality of Clark-type quotient Hardy modules, we need a preliminary lemma.
\begin{lem}\label{lem:MobiusTransform}
    Let $\mathcal{N}$ be a quotient module of $H^2(\mathbb{D}^d)$, and $\xi$ be an inner function on $\mathbb{D}^d$ such that $S_\xi$ is $(1,\infty)$-essentially unitary, and $\varphi_a:z\mapsto\frac{a-z}{1-\bar{a}z}$ be a M\"obius transform. Then $S_{\varphi_a(\xi)}$ is also $(1,\infty)$-essentially unitary.
\end{lem}
\begin{proof}
    Because $S_\xi$ is essentially unitary, it holds that $\sigma_e(S_\xi)\subset\mathbb{T}$. By the Spectral Mapping Theorem, the essential spectrum $$\sigma_e(S_{\varphi_a(\xi)})=\varphi_a\big(\sigma_e(S_\xi)\big)\subset\mathbb{T},$$ which implies the Fredholmness of $S_{\varphi_a(\xi)}$.

    Since $S_\xi$ is $(1,\infty)$-essentially unitary, we have $$P_{\mathcal{N}}M_\xi^*P_{\mathcal{N}^\perp}M_\xi P_{\mathcal{N}}=(\Id_{\mathcal{N}}-S_\xi^*S_\xi)P_{\mathcal{N}}\in\mathcal{L}^{(1,\infty)},$$
    which gives $P_{\mathcal{N}^\perp}M_\xi P_{\mathcal{N}}\in\mathcal{L}^{(2,\infty)}$. We set $c=\|P_{\mathcal{N}^\perp}M_\xi P_{\mathcal{N}}\|_{(2,\infty)}$. Observe that
    \begin{eqnarray*}
        \|P_{\mathcal{N}^\perp}M_\xi^{k+1}P_{\mathcal{N}}\|_{(2,\infty)}&=
        &\|P_{\mathcal{N}^\perp}M_\xi^k(P_{\mathcal{N}^\perp}+P_{\mathcal{N}})M_\xi P_{\mathcal{N}}\|_{(2,\infty)}\\
        &\leq&\|P_{\mathcal{N}^\perp}M_\xi^kP_{\mathcal{N}^\perp}M_\xi P_{\mathcal{N}}\|_{(2,\infty)}
        +\|P_{\mathcal{N}^\perp}M_\xi^kP_{\mathcal{N}}M_\xi P_{\mathcal{N}}\|_{(2,\infty)}\\
        &\leq&c+\|P_{\mathcal{N}^\perp}M_\xi^kP_{\mathcal{N}}\|_{(2,\infty)},
    \end{eqnarray*}
    then by induction we get
    $$\|P_{\mathcal{N}^\perp}M_\xi^kP_{\mathcal{N}}\|_{(2,\infty)}\leq kc,~\forall k\in\mathbb{N}.$$
    Consequently,
    \begin{eqnarray*}
        \|P_{\mathcal{N}^\perp}M_{\varphi_a(\xi)}P_{\mathcal{N}}\|_{(2,\infty)}&=&\|P_{\mathcal{N}^\perp}
        (a-M_\xi)(1-\bar{a}M_\xi)^{-1}P_{\mathcal{N}}\|_{(2,\infty)}\\
        &=&\|\sum_{k=0}^\infty a\bar{a}^kP_{\mathcal{N}^\perp}M_\xi^kP_{\mathcal{N}}-\sum_{k=0}^\infty\bar{a}^kP_{\mathcal{N}^\perp}
        M_\xi^{k+1}P_{\mathcal{N}}\|_{(2,\infty)}\\
        &\leq&\sum_{k=0}^\infty|a|^{k+1}\|P_{\mathcal{N}^\perp}M_\xi^kP_{\mathcal{N}}\|_{(2,\infty)}
        +\sum_{k=0}^\infty|a|^k\|P_{\mathcal{N}^\perp}M_\xi^{k+1}P_{\mathcal{N}}\|_{(2,\infty)}\\
        &\leq&\sum_{k=0}^\infty kc|a|^{k+1}+\sum_{k=0}^\infty(k+1)c|a|^k\\
        &<&\infty,
    \end{eqnarray*}
    which implies $P_{\mathcal{N}^\perp}M_{\varphi_a(\xi)}P_{\mathcal{N}}\in\mathcal{L}^{(2,\infty)}$. Moreover
    \begin{equation}\label{eq:eunitary3}
        \Id_{\mathcal{N}}-S_{\varphi_a(\xi)}^*S_{\varphi_a(\xi)}=P_{\mathcal{N}}M_{\varphi_a(\xi)}^*
        P_{\mathcal{N}^\perp}M_{\varphi_a(\xi)}P_{\mathcal{N}}\mid_{\mathcal{N}}\in\mathcal{L}^{(1,\infty)}.
    \end{equation}
    Let $S_{\varphi_a(\xi)}=W|S_{\varphi_a(\xi)}|$ be the polar decomposition, where $|S_{\varphi_a(\xi)}|=\big(S_{\varphi_a(\xi)}^*S_{\varphi_a(\xi)}\big)^{1/2}$ and $W$ is the partial isometry that maps $\big(\ker S_{\varphi_a(\xi)}\big)^\perp$ unitarily onto $\mathrm{ran}~S_{\varphi_a(\xi)}$. Since $S_{\varphi_a(\xi)}$ is Fredholm, the subspace $\range(\Id_{\mathcal{N}}-WW^*)=\coker S_{\varphi_a(\xi)}$ is finite dimensional. Then it follows that
    $$\Id_{\mathcal{N}}-S_{\varphi_a(\xi)}S_{\varphi_a(\xi)}^*=W(\Id_{\mathcal{N}}-S_{\varphi_a(\xi)}^*S_{\varphi_a(\xi)})
    W^*+(\Id_{\mathcal{N}}-WW^*)\in\mathcal{L}^{(1,\infty)},$$
    which completes the proof of the lemma.
\end{proof}
For a tuple $z\in\mathbb{C}^d$ and a multi-index $\alpha\in\mathbb{N}^d$, we write $z^\alpha=\prod_{i=1}^dz_i^{\alpha_i}$ for brevity. Similarly, $B^\alpha(z)=\prod_{i=1}^dB_i^{\alpha_i}(z_i)$. To prove the essential normality of Clark-type quotient modules, we need a deeper interpretation of \cite[Proposition 3.1]{WZ3}.
\begin{prop}\label{prop:eu_quasi-homogeneous_distinguished}
    Let $I$ be a $\mathbf{K}$-quasi-homogeneous distinguished ideal of $\poly$, then for $1\leq i\leq d$, $S_{z_i,[I]^\perp}$ is $(1,\infty)$-essentially unitary, and $P_{[I]}M_{z_i}P_{[I]^\perp}\in\mathcal{L}^{(2,\infty)}.$
\end{prop}
\begin{proof}
    By \cite[Proposition 3.1]{WZ1} we have $\dim_{\mathbb{C}}Z(I)\leq1$. In the case $I$ is primary and homogeneous, $Z(I)\cap\mathbb{D}^d$ is a disc, hence \cite[Theorem 2.17]{WZ1} shows that $P_{[I]}M_{z_i}P_{[I]^\perp}\in\mathcal{L}^{(2,\infty)}$.

    In the case $I$ is homogeneous, let $I=\bigcap_{k=1}^n I_k$ be an irredundant primary decomposition of $I$. Then $P_{[I_k]}M_{z_i}P_{[I_k]^\perp}\in\mathcal{L}^{(2,\infty)}$ for $1\leq k\leq n$.
    According to the proof of \cite[Theorem 3.2]{WZ1}, there is a positive constant $c$ such that $P_{[I]^\perp}\leq c\sum_{k=1}^nP_{[I_k]^\perp}$, hence there is a bounded linear operator $A$ such that $P_{[I]^\perp}=\sum_{k=1}^nP_{[I_k]^\perp}A$. Consequently
    $$P_{[I]}M_{z_i}P_{[I]^\perp}=P_{[I]}M_{z_i}\sum_{k=1}^nP_{[I_k]^\perp}A=P_{[I]}\sum_{k=1}^nP_{[I_k]}M_{z_i}P_{[I_k]^\perp}A\in\mathcal{L}^{(2,\infty)}.$$

    In the general case, the linear mapping $\Psi_{\mathbf{K}}:H^2(\mathbb{D}^d)\to H^2(\mathbb{D}^d)$ defined by $z^\alpha\mapsto z^{(\alpha_1K_1,\ldots,\alpha_dK_d)}$
    is an isometry with $\Psi_{\mathbf{K}}M_{z_i}=M_{z_i}^{K_i}\Psi_{\mathbf{K}},~1\leq i\leq d$. Clearly the ideal $\tilde{I}$ generated by $\Psi_{\mathbf{K}}(I)$ is homogeneous and distinguished. The argument in the preceding paragraph shows that $P_{[\tilde{I}]}M_{z_i}P_{[\tilde{I}]^\perp}$ lies in $\mathcal{L}^{(2,\infty)}$, and recursively we can prove that
    $$P_{[\tilde{I}]}M_{z_i}^nP_{[\tilde{I}]^\perp}=P_{[\tilde{I}]}M_{z_i}^{n-1}P_{[\tilde{I}]}M_{z_i}P_{[\tilde{I}]^\perp}+P_{[\tilde{I}]}M_{z_i}^{n-1}P_{[\tilde{I}]^\perp}M_{z_i}P_{[\tilde{I}]^\perp}\in\mathcal{L}^{(2,\infty)}$$
    for $n=2,3,\ldots$. In particular $P_{[\tilde{I}]}M_{z_i}^{K_i}P_{[\tilde{I}]^\perp}\in\mathcal{L}^{(2,\infty)}$, which gives $P_{\Psi_{\mathbf{K}}[I]}M_{z_i}^{K_i}P_{\Psi_{\mathbf{K}}[I]^\perp}\in\mathcal{L}^{(2,\infty)}$ i.e. $P_{[I]}M_{z_i}P_{[I]^\perp}\in\mathcal{L}^{(2,\infty)}.$

    From the following facts
    $$\Id_{[I]^\perp}-S_{z_i,[I]^\perp}^*S_{z_i,[I]^\perp}=P_{[I]^\perp}M_{z_i}^*P_{[I]}M_{z_i}P_{[I]^\perp}\in\mathcal{L}^{(1,\infty)}$$
    we see that $S_{z_i,[I]^\perp}$ is $(1,\infty)$-essentially isometry, and then the essential normality of $S_{z_i,[I]^\perp}$ \cite[Theorem 3.1]{WZ3} ensures its Fredholmness. Let $S_{z_i}=W|S_{z_i}|$ be the polar decomposition of $S_{z_i}$, where $W$ is the partial isometry mapping $(\ker S_{z_i})^\perp$ onto $\range S_{z_i}$. Then the fact
    $$\Id_{[I]^\perp}-S_{z_i,[I]^\perp}S_{z_i,[I]^\perp}^*=W(\Id_{[I]^\perp}-S_{z_i,[I]^\perp}^*S_{z_i,[I]^\perp})W^*+P_{\ker S_{z_i}^*}\in\mathcal{L}^{(1,\infty)}$$
    implies that $S_{z_i,[I]^\perp}$ is $(1,\infty)$-essentially co-isometry, completing the proof of the proposition.
\end{proof}
The following result gives the essential normality of a large class of inhomogeneous quotient Hardy modules.
\begin{thm}\label{thm:enormal_Hardy}
    Let $I\subset\poly$ be a Clark-type ideal, then the quotient module $[I]^\perp$ of $H^2(\mathbb{D}^d)$ is $(1,\infty)$-essentially normal. Moreover, $S_{z_i}$ is $(1,\infty)$-essentially unitary for $i=1,\ldots,d$.
\end{thm}
\begin{proof}
    Let $J,B_1,\ldots,B_d$ be as in Definition \ref{def:Clarktype}. Since each $B_i$ is a nonconstant finite Blaschke product, the multiplication operator $M_{B_i}$ on $H^2(\mathbb{D})$ is an isometry with $\dim_{\mathbb{C}}\ker M_{B_i}^*<\infty$. Write $H_0=\bigotimes_{i=1}^d\ker M_{B_i}^*$ which is finite dimensional, then by the Wold decomposition for isometries we have
    $$H^2(\mathbb{D}^d)=\bigoplus_{\alpha\in\mathbb{N}^d}M_{B^\alpha}H_0.$$
    For $p\in\poly$ and $f\in H^2(\mathbb{D}^d)$ we write $A_pf=p(B)f$, then $A$ defines a $\poly$-module structure on $H^2(\mathbb{D}^d)$. Define
    \begin{eqnarray*}
        U:H^2(\mathbb{D}^d)\otimes H_0&\to&H^2(\mathbb{D}^d)\\
        \sum_{\alpha\in\mathbb{N}^d}z^\alpha\otimes f_\alpha&\mapsto&\sum_{\alpha\in\mathbb{N}^d}B^\alpha f_\alpha,
    \end{eqnarray*}
    then $U$ is unitary and $U(M_p\otimes\Id_{H_0})=A_pU$. For $f\in H^2(\mathbb{D}^d)\otimes H_0$, we have
    \begin{eqnarray*}
        Uf\in[I]^\perp&\Leftrightarrow&M_p^*Uf=0,~\forall p\in I\\
        &\Leftrightarrow&M_{q(B)}^*Uf=0,~\forall q\in J\\
        &\Leftrightarrow&A_q^*Uf=0,~\forall q\in J\\
        &\Leftrightarrow&U(M_q\otimes\Id_{H_0})^*f=0,~\forall q\in J\\
        &\Leftrightarrow&(M_q\otimes\Id_{H_0})^*f=0,~\forall q\in J\\
        &\Leftrightarrow&f\in[J]^\perp\otimes H_0,
    \end{eqnarray*}
    and therefore $[I]^\perp=U([J]^\perp\otimes H_0)$. In \cite{WZ1,WZ3} it is proven that $[J]^\perp$ is essentially normal, and Proposition \ref{prop:eu_quasi-homogeneous_distinguished} gives
    $$P_{[J]}M_{z_i}P_{[J]^\perp}\in\mathcal{L}^{(2,\infty)},$$
    and
    $$P_{[J]^\perp}-P_{[J]^\perp}M_{z_i}M_{z_i}^*P_{[J]^\perp}\in\mathcal{L}^{(1,\infty)}.$$
    Write $B_i(z_i)=\varphi_a(z_i)\tilde{B}_i(z_i)$ where $\varphi_a$ is a M\"obius transform and $\tilde{B}_i$ is a Blaschke product. It follows that
    \begin{equation}\label{eq:conner_operator}
        P_{[I]}M_{B_i(z_i)}P_{[I]^\perp}=U\big((P_{[J]}M_{z_i}P_{[J]^\perp})\otimes\Id_{H_0}\big)U^*\in\mathcal{L}^{(2,\infty)},
    \end{equation}
    and
    \begin{eqnarray}\label{eq:eunitary1}
        P_{[I]^\perp}-S_{\varphi_a(z_i)}S_{\varphi_a(z_i)}^*&=&P_{[I]^\perp}-P_{[I]^\perp}M_{\varphi_a(z_i)}M_{\varphi_a(z_i)}^*P_{[I]^\perp}\notag\\
        &=&P_{[I]^\perp}P_{\ker M_{\varphi_a(z_i)}^*}P_{[I]^\perp}\notag\\
        &\leq&P_{[I]^\perp}P_{\ker M_{B_i(z_i)}^*}P_{[I]^\perp}\notag\\
        &=&P_{[I]^\perp}-P_{[I]^\perp}M_{B_i(z_i)}M_{B_i(z_i)}^*P_{[I]^\perp}\notag\\
        &=&U\big((P_{[J]^\perp}-P_{[J]^\perp}M_{z_i}M_{z_i}^*P_{[J]^\perp})\otimes\Id_{H_0}\big)U^*\notag\\
        &\in&\mathcal{L}^{(1,\infty)}.
    \end{eqnarray}
    Similarly, from (\ref{eq:conner_operator}) and
    \begin{eqnarray*}
        &&P_{[I]^\perp}M_{B_i(z_i)}^*P_{[I]}M_{B_i(z_i)}P_{[I]^\perp}\\
        &=&P_{[I]^\perp}-P_{[I]^\perp}M_{B_i(z_i)}^*P_{[I]^\perp}M_{B_i(z_i)}P_{[I]^\perp}\\
        &=&P_{[I]^\perp}-P_{[I]^\perp}M_{\varphi_a(z_i)}^*P_{[I]^\perp}M_{\varphi_a(z_i)}P_{[I]^\perp}+P_{[I]^\perp}M_{\varphi_a(z_i)}^*(P_{[I]^\perp}-M_{\tilde{B}_1(z)}^*P_{[I]^\perp}M_{\tilde{B}_1(z)})M_{\varphi_a(z_i)}P_{[I]^\perp}\\
        &\geq&P_{[I]^\perp}-P_{[I]^\perp}M_{\varphi_a(z_i)}^*P_{[I]^\perp}M_{\varphi_a(z_i)}P_{[I]^\perp}\\
        &=&P_{[I]^\perp}M_{\varphi_a(z_i)}^*P_{[I]}M_{\varphi_a(z_i)}P_{[I]^\perp},
    \end{eqnarray*}
    we obtain
    \begin{equation}\label{eq:eunitary2}
        \Id_{[I]^\perp}-S_{\varphi_a(z_i)}^*S_{\varphi_a(z_i)}=P_{[I]^\perp}M_{\varphi_a(z_i)}^*P_{[I]}M_{\varphi_a(z_i)}P_{[I]^\perp}\mid_{[I]^\perp}\in\mathcal{L}^{(1,\infty)}.
    \end{equation}
    Combining (\ref{eq:eunitary1}) and (\ref{eq:eunitary2}), we conclude that $S_{\varphi_a(z_i)}$ is $(1,\infty)$-essentially unitary. Then it follows from Lemma \ref{lem:MobiusTransform} that $S_{z_i}=S_{\varphi_{a}\circ\varphi_a(z_i)}$ is $(1,\infty)$-essentially unitary.
\end{proof}
In the proof of Theorem \ref{thm:enormal_Hardy}, if one of the inner functions $B_1,\ldots,B_d$, for instance $B_1$, is not a finite Blaschke product, then
$$H_0=\ker M_{B_1}^*\otimes\bigotimes_{i=2}^d\ker M_{B_i}^*$$
is an infinite dimensional subspace. The unitary $U$ builds up the equivalence between $S_{B_2(z_2),[I]^\perp}$ and $S_{z_2,[J]^\perp}\otimes\Id_{H_0}$, which are not essentially normal. Consequently, $S_{z_2,[I]^\perp}$ cannot be essentially normal. This leads to the following result.
\begin{cor}\label{cor:ess}
    Let $\eta=\big(\eta_1(z_1),\ldots,\eta_d(z_d)\big)$ be a $d$-tuple of nonconstant inner functions, and $J\subset\poly$ be a distinguished quasi-homogeneous ideal such that $\dim Z(J)\cap\mathbb D^d>0$. Then $[J_\eta]^\perp$ is essentially normal if and only if all the $\eta_i$'s are finite Blaschke products.
\end{cor}
Using techniques similar to that in \cite[Section 5]{WZ3}, Theorem \ref{thm:enormal_Hardy} can be generalized to the case of integer-weighted Bergman modules. Recall that $[I]_\alpha^\perp$ denotes the quotient module of
$A_{\alpha}(\mathbb{D}^2)$ determined by the ideal $I$, and $S_{f,[I]_\alpha^\perp}=P_{[I]_\alpha^\perp}M_f|_{[I]_\alpha^\perp}.$
\begin{cor}\label{cor:enormal}
    Let $I\subset\mathbb{C}[z,w]$ be a Clark-type ideal, then for every integer $\alpha\geq-1$, the quotient weighted Bergman module $[I]_\alpha^\perp$ is $(1,\infty)$-essentially normal. Moreover, both $S_{z,[I]_\alpha^\perp}$ and $S_{w,[I]_\alpha^\perp}$ are $(1,\infty)$-essentially unitary.
\end{cor}
\begin{proof}
    The case $\alpha=-1$ is included in Theorem \ref{thm:enormal_Hardy}, so in the sequel we assume $\alpha\geq0$. According to \cite[Proposition 2.5]{GWZ}, it suffices to prove that both $S_{z,[I]^\perp}$ and $S_{w,[I]^\perp}$ are $(1,\infty)$-essential co-isometries.

    By Definition \ref{def:Clarktype}, there is a distinguished quasi-homogeneous ideal $J\in\mathbb{C}[z,w]$ and nonconstant finite Blaschke products $B_1,B_2$, such that $I$ is generated by the numerators of functions in $J_B$. Let $J_1$ be the quasi-homogeneous ideal of $\mathbb{C}[z_1,\ldots,z_{\alpha+2},w_1,\ldots,w_{\alpha+2}]$ generated by
    $$\{q(z_i,w_i),~z_i-z_j,~w_i-w_j:q\in J,~1\leq i,j\leq\alpha+2\},$$
    and $I_1\subset\mathbb{C}[z_1,\ldots,z_{\alpha+2},w_1,\ldots,w_{\alpha+2}]$ be the ideal generated by the numerators of functions in $$\{q\big(B_1(z_1),\ldots,B_1(z_{\alpha+2}),B_2(w_1),\ldots,B_2(w_{\alpha+2})\big):q\in J_1\}.$$
    By Theorem \ref{thm:enormal_Hardy}, on the quotient Hardy module $[I_1]^\perp\subset H^2(\mathbb{D}^{2\alpha+4})$, the operators $S_{z_i,[I_1]^\perp}$ and $S_{w_i,[I_1]^\perp}$ are $(1,\infty)$-essentially unitary.

    Let $I_2\subset\mathbb{C}[z_1,\ldots,z_{\alpha+2},w_1,\ldots,w_{\alpha+2}]$ be the ideal generated by the numerators of functions in
    $$\{p(z_i,w_i),~z_i-z_j,~w_i-w_j:1\leq i,j\leq\alpha+2,~p\in I\}.$$
    It is routine to verify that $B_1(z_i)-B_1(z_j)$ and $B_2(w_i)-B_2(w_j)$ lie in $[I_2]$ for $1\leq i,j\leq\alpha+2$, and $q\big(B_1(z_i),B_2(w_i)\big)\in[I_2]$ for $q\in J$.

    Hence $[I_1]\subset[I_2]$ or equivalently $[I_2]^\perp\subset[I_1]^\perp$. Consequently, both $S_{z_i,[I_2]^\perp}$ and $S_{w_i,[I_2]^\perp}$ are $(1,\infty)$-essentially co-isometric. By \cite{Cla,WZ1}, for the ideal $I_3$ generated by $\{z_i-z_j,~w_i-w_j:1\leq i\leq\alpha+2\}$, the quotient Hardy module $[I_3]^\perp$ is unitarily equivalent to the weighted Bergman module $\mathcal{A}_\alpha(\mathbb{D}^2)$, via the mapping
    $$f(z_1,\ldots,z_{\alpha+2},w_1,\ldots,w_{\alpha+2})\mapsto f(z,\ldots,z,w,\ldots,w).$$
    From the following facts
    \begin{eqnarray*}
        f\in[I_2]^\perp&\Leftrightarrow&M_{p(z_i,w_i)}^*f=M_{z_i-z_j}^*f=M_{w_i-w_j}^*f=0,~1\leq i,j\leq\alpha+2,~\forall p\in I\\
        &\Leftrightarrow&f\in[I_3]^\perp,~M_{p(z,w)}^*f(z,\ldots,z,w,\ldots,w)=0,~\forall p\in I\\
        &\Leftrightarrow&f\in[I_3]^\perp,~f(z,\ldots,z,w,\ldots,w)\in[I]_\alpha^\perp,
    \end{eqnarray*}
    we find that $[I_2]^\perp$ as a quotient module is unitarily equivalent to $[I]_\alpha^\perp$. Since $S_{z_i,[I_2]^\perp}$ and $S_{w_i,[I_2]^\perp}$ are $(1,\infty)$-essentially co-isometric, then so are $S_{z,[I]_\alpha^\perp}$ and $S_{w,[I]_\alpha^\perp}$, completing the proof of the corollary.
\end{proof}
The following theorem omits the requirement to the algebraic structure of Clark-type ideal in Corollary \ref{cor:enormal}, and relates the essential normality of quotient modules to the varieties.
\begin{thm}\label{thm:enormal}
    If $I\subset\mathbb{C}[z,w]$ is a distinguished ideal with Clark-type variety, then for every integer $\alpha\geq-1$, the quotient weighted Bergman module $[I]_\alpha^\perp$ is $(1,\infty)$-essentially normal. Moreover, $S_{z,[I]_\alpha^\perp}$ and $S_{w,[I]_\alpha^\perp}$ are $(1,\infty)$-essentially unitary.
\end{thm}
\begin{proof}Without loss of generality, assume that any algebraic branches of $Z(I)$ intersects with $\mathbb D^2$.
    According to \cite[Proposition 2.5]{GWZ}, it suffices to prove that both $S_{z,[q]^\perp}$ and $S_{w,[q]^\perp}$ are $(1,\infty)$-essentially co-isometry.

    By Definition \ref{def:Clarktype}, there exists a radical distinguished quasi-homogeneous ideal $J=p\mathbb{C}[z,w]$ and finite Blaschke products $B_1,B_2$ such that
    $$Z(I)\cap\mathbb{D}^2=Z(J_B)\cap\mathbb{D}^2=Z(p_B)\cap\mathbb{D}^2.$$
    Write $p_B=q/r$ where $q$ and $r$ are coprime polynomials, then $q$ is distinguished.

    The Hilbert's Nullstellensatz ensures the existence of an integer $N$ such that $q^N\in I$. Clearly
    \begin{eqnarray*}
        &&\{f\in\mathbb{C}[z,w]:f=gp_B^N,~g\in\mathbb{C}[z,w]\}\\
        &=&\{f\in\mathbb{C}[z,w]:f=gq^N/r^N,~g\in\mathbb{C}[z,w]\}\\
        &=&(q^N),
    \end{eqnarray*}
    Consequently $(q^N)$ a Clark-type ideal. Corollary \ref{cor:enormal} shows that $S_{z,[q^N]_\alpha^\perp}$ is $(1,\infty)$-essentially unitary, and therefore
    $$\Id_{[I]_\alpha^\perp}-S_{z,[I]_{\alpha}^\perp}S_{z,[I]_{\alpha}^\perp}^*=
    P_{[I]_\alpha^\perp}\left(\Id_{[q^N]_\alpha^\perp}-S_{z,[q^N]_\alpha^\perp}S_{z,[q^N]_\alpha^\perp}^*\right)\mid_{[I]_\alpha^\perp}\in\mathcal{L}^{(1,\infty)},$$
    and symmetrically
    $$\Id_{[I]_\alpha^\perp}-S_{w,[I]_{\alpha}^\perp}S_{w,[I]_{\alpha}^\perp}^*\in\mathcal{L}^{(1,\infty)}$$
    as desired.
\end{proof}
By \cite[Corollary 6.4]{GWZ}, if $I_1$ and $I_2$ are distinguished ideals of $\mathbb{C}[z,w]$ such that $Z(I_1)$ and $Z(I_2)$ do not intersect on $\mathbb{T}^2$, and both $[I_1]_\alpha^\perp$ and $[I_2]_\alpha^\perp$ are essentially normal, then $[I_1\cap I_2]_\alpha^\perp$ is also essentially normal. Combining this fact with Theorem \ref{thm:enormal}, we immediately get the following result.
\begin{cor}
    Let $I\subset\mathbb{C}[z,w]$ be an ideal which is the intersection of ideals $I_1,\ldots,I_n$ with Clark-type varieties, such that $Z(I_i)\cap Z(I_j)\cap\mathbb{T}^2=\emptyset$ whenever $i\neq j$. Then for each integer $\alpha>-1$, the quotient weighted Bergman module $[I]_\alpha^\perp$ is essentially normal.
\end{cor}
\begin{rem}\label{rem:primary_decomposition}
    For a distinguished ideal $I$ of $\poly$, since each component of $Z(I)$ intersects with $\mathbb{D}^d$, in each irredundant primary decomposition $I=\bigcap I_k$, either $Z(I_k)$ is a single point in $\mathbb{D}^d$ or $\dim_{\mathbb{C}}Z(I_k)\geq1$. In the former case $I_k$ is finite codimensional in $\poly$, and in the latter case $I_k$ is a distinguished ideal. By the theory of primary decomposition, the primary components $I_k$ with $\dim_{\mathbb{C}}Z(I_k)\geq1$ are uniquely determined by $I$. Because the essential normality of $[I]^\perp$ is not impacted by the primary components with finite codimension, for convenience we only refer to their intersection, and write the primary decomposition as $I=\bigcap_{k=0}^mI_k$, where $I_0$ is a finite codimensional ideal, and $I_1,\ldots,I_m$ are primary ideals with $\dim_{\mathbb{C}}Z(I_k)\geq1$.
\end{rem}
\section{K-homology for Clark-type quotient modules}

~~~~In this section, we consider the $K$-homology defined by Clark-type quotient modules. By the BDF-theory \cite{BDF}, when a quotient module $[I]^\perp$ is essentially normal, we have the following short exact sequence
\begin{eqnarray}\label{exact-sequence}
0\to {\cal K}\to C^*([I]^\perp)+{\cal K}\to C(\sigma_e(S_{z_1},\cdots,S_{z_d}))\to 0,
\end{eqnarray}
where $C^*([I]^\perp)$ is the $C^*$-algebra generated by $\{\Id_{[I]^\perp}, S_{z_1},\cdots, S_{z_d}\}$, and $\sigma_e(S_{z_1},\cdots,S_{z_d})$ is the Taylor joint essential spectrum \cite{Cur,Tay1} of the tuple $(S_{z_1},\cdots,S_{z_d})$. We will show that the short exact sequence (\ref{exact-sequence}) is not splitting, i.e. $C^*([I]^\perp)+{\cal K}$ is not the direct product of $\mathcal{K}$ and $C(\sigma_e(S_{z_1},\cdots,S_{z_d}))$, then (\ref{exact-sequence}) gives a nontrivial element $e_I$ in the K-homology group $K\big(\sigma_e(S_{z_1},\cdots,S_{z_d})\big)$. To this aim, we first show that
$$\sigma_e(S_{z_1},\cdots,S_{z_d})=Z(I)\cap\partial\mathbb{D}^d$$
for Clark-type quotient modules $[I]^\perp$.

For $z\in\mathbb{C}^d$ and a subset $S$ of $\{1,\ldots,d\}$, we write $z_S=\prod_{i\in S}z_i$. For a Hilbert $\poly$-module $\mathcal{H}$, we set
\begin{eqnarray}\label{Def:Core-Operator}C_\mathcal{H}=\sum_{S\subset\{1,\ldots,d\}}(-1)^{\card S}M_{z_S}M_{z_S}^*,\end{eqnarray}
then for a submodule $\mathcal{M}$ of $H^2(\mathbb{D}^d)$, $C_{\mathcal{M}}$ is just the core operator defined by Guo and Yang \cite{GY}. It is routine to verify that $C_{H^2(\mathbb{D}^d)}=1\otimes1$ is of rank $1$. For the weighted Bergman module $\mathcal{A}_\alpha(\mathbb{D}^d)$ where $\alpha\in(-1,\infty)$, an orthonormal basis is
$$\left\{\sqrt{\prod_{i=1}^d\frac{\Gamma(\beta_i+\alpha+2)}{\Gamma(\beta_i+1)\Gamma(\alpha+2)}}z^\beta:\beta\in\mathbb{N}^d\right\}.$$
Direct computations verify that
$$C_{\mathcal{A}_\alpha(\mathbb{D}^d)}z^\beta=\frac{(\alpha+1)^d}{\prod_{i=1}^d(\beta_i+\alpha+1)}z^\beta,~\forall\beta\in\mathbb{N}^d,$$
which ensures the compactness of $C_{\mathcal{A}_\alpha(\mathbb{D}^d)}$.
\begin{prop}\label{prop:espectrum}
    Let $I\subset\poly$ be an ideal such that
    \begin{equation}\label{eq:varietyboundary}
        \partial\big(Z(I)\cap\mathbb{D}^d\big)=Z(I)\cap\partial\mathbb{D}^d,
    \end{equation}
    and $\alpha\geq-1$ be a real number, such that the quotient module $[I]_\alpha^\perp$ of $\mathcal{A}_\alpha(\mathbb{D}^d)$ is essentially normal, then the Taylor essential spectrum \cite{Tay1,Tay2}
    $$\sigma_e(S_{z_1},\ldots,S_{z_d})=Z(I)\cap\partial\mathbb{D}^d.$$
\end{prop}
\begin{proof}
    Since $p(S_{z_1},\ldots,S_{z_d})=0$ for all polynomials $p\in I$, the Spectral Mapping Theorem \cite{Tay2} gives
    $$\sigma_e(S_{z_1},\ldots,S_{z_d})\subset\sigma(S_{z_1},\ldots,S_{z_d})\subset Z(I).$$
    From
   \begin{equation}
        C_{[I]_\alpha^\perp}=P_{[I]_\alpha^\perp}C_{\mathcal{A}_\alpha(\mathbb{D}^d)}\mid_{[I]_\alpha^\perp}
    \end{equation}
    we obtain the compactness of $C_{[I]_\alpha^\perp}$. Then since $C^*([I]_\alpha^\perp)$ is essentially commutative, it follows from the definition of $C_{[I]_\alpha^\perp}$ that $\sigma_e(S_{z_1},\ldots,S_{z_d})\subset\partial\mathbb{D}^d$. Therefore we have proved
    $$\sigma_e(S_{z_1},\ldots,S_{z_d})\subset Z(I)\cap\partial\mathbb{D}^d.$$

    For the inverse inclusion, let $\lambda\in Z(I)\cap\partial\mathbb{D}^d$. By \cite[Corollary 3.9]{Cur}, it is sufficient to show that $T=\sum\limits_{i=1}^{d}(\lambda_i-S_{z_i})(\lambda_i-S_{z_i})^*$ is not Fredholm. If $T$ is Fredholm, then by the positivity of $T$, we can find an invertible positive operator $B$ and a compact $K$ such that $T=B+K$. Take a sequence $\{\mu^{(n)}\}$ in $Z(I)\cap\mathbb{D}^d$ that converges to $\lambda$. Let $k_{\mu^{(n)}}=\frac{K_{\mu^{(n)}}}{\|K_{\mu^{(n)}}\|}$ be the normalized reproducing kernel, then $\{k_{\mu^{(n)}}\}$ converges to $0$ weakly. Hence there is a positive number $c$ such that
    \begin{eqnarray*}
    \liminf_{n\to\infty}\langle T k_{\mu^{(n)}}, k_{\mu^{(n)}}\rangle=\liminf_{n\to\infty}\langle(B+K)k_{\mu^{(n)}}, k_{\mu^{(n)}}\rangle=\liminf_{n\to\infty}\langle Bk_{\mu^{(n)}}, k_{\mu^{(n)}}\rangle\geq c.
    \end{eqnarray*}
    However, since $\mu^{(n)}\in Z(I)\cap\mathbb D^2$, it holds that $k_{\mu^{(n)}}\in [I]^\perp$, and consequently
\begin{eqnarray}
\lim\limits_{n\to\infty}\langle T k_{\mu^{(n)}},k_{\mu^{(n)}}\rangle=\lim\limits_{n\to\infty} |\lambda-\mu^{(n)}|^2=0,
\end{eqnarray}
contradicting to the previous inequality. Therefore $T$ is not Fredholm and $\lambda\in\sigma_e(S_{z_1},\ldots,S_{z_d})$, completing the proof of the proposition.
\end{proof}
For a quasi-homogeneous distinguished ideal $J\subset\poly$, \cite{WZ3} proves the essential normality of $[J]^\perp$. By quasi-homogeneity we have
$$\partial\big(Z(J)\cap\mathbb{D}^d\big)=Z(J)\cap\mathbb{T}^d=Z(J)\cap\partial\mathbb{D}^d.$$
For the Clark-type variety $V=Z(J_B)\cap\mathbb{D}^d$ where $B(z)=\big(B_1(z_1,\ldots,B_d(z_d)\big)$ is a tuple of finite Blaschke product, clearly we have $\partial V\subset Z(J_B)\cap\partial\mathbb{D}^d$. Conversely, for each $z\in Z(J_B)\cap\partial\mathbb{D}^d$, it holds that $p\big(B(z)\big)=0,~\forall p\in J$ while $B(z)\in\partial\mathbb{D}^d$, i.e.
$$B(z)\in Z(J)\cap\partial\mathbb{D}^d=\partial\big(Z(J)\cap\mathbb{D}^d\big)=Z(J)\cap\mathbb{T}^d.$$
Consequently $z\in\mathbb{T}^d$, and there is a sequence $\{\lambda^{(n)}\}$ in $Z(J)\cap\mathbb{D}^d$ that converges to $B(z)$. That is, for each $1\leq i\leq d$, it holds that $\lim_{n\to\infty}\lambda_i^{(n)}=B_i(z_i)$. Since $B_i$ has no singularities on $\mathbb{T}$, we can choose $\mu_i^{(n)}\in B_i^{-1}(\lambda_i^{(n)})$ such that $\lim_{n\to\infty}\mu_i^{(n)}=z_i$. Write $\mu^{(n)}=(\mu_1^{(n)},\ldots,\mu_d^{(n)})$, then $B(\mu^{(n)})=\lambda^{(n)}$ and a priori $\mu^{(n)}\in Z(J_B)\cap\mathbb{D}^d=V$, and $\lim_{n\to\infty}\mu^{(n)}=z$. This shows that $Z(J_B)\cap\partial\mathbb{D}^d\subset\partial V$.

Therefore the condition (\ref{eq:varietyboundary}) in Proposition \ref{prop:espectrum} is satisfied by Clark-type ideals. This together with Theorem \ref{thm:enormal_Hardy} and Theorem \ref{thm:enormal} derives the following corollaries on the Taylor essential spectrums of the quotient modules.
\begin{cor}
    Let $I\subset\poly$ be a Clark-type ideal, then
    $$\sigma_e(S_{z_1},\ldots,S_{z_d})=Z(I)\cap\partial\mathbb{D}^d.$$
\end{cor}
As a special case for $d=2,$ the following example.
\begin{exam}
    Let $I\subset\mathbb{C}[z,w]$ be an ideal with Clark-type variety, and $\alpha\geq-1$ be an integer, then
    $$\sigma_e(S_{z,[I]_\alpha^\perp},S_{w,[I]_\alpha^\perp})=Z(I)\cap\partial\mathbb{D}^2.$$
\end{exam}
By Theorem \ref{thm:enormal} and \cite[Theorem 5.2]{GWZ}, we can prove the non-triviality of the $C^*$-extensions associated to the quotient modules with Clark-type variety.
\begin{thm}
    Let $I\subset\mathbb{C}[z,w]$ be an ideal with Clark-type variety, and $\alpha\geq-1$ be an integer, then the short exact sequence
    $$0\to \cal{K} \hookrightarrow C^*({[I]_\alpha^\perp})+\mathcal{K}\to C(Z(I)\cap \partial\mathbb D^2)\to 0$$
    is not splitting.
\end{thm}
To prove the non-triviality of the $C^*$-extensions corresponding to the Clark-type quotient modules in the high dimensional cases, the following technical lemma is needed.
\begin{lem}\label{lem:nonsplit}
    Let $J\subset\poly$ be a distinguished $\mathbf{K}$-quasi-homogeneous ideal, then $S_{z_1,[J]^\perp}$ is a Fredholm operator with nonzero index.
\end{lem}
\begin{proof}
    Proposition \ref{prop:eu_quasi-homogeneous_distinguished} shows that $S_{z_i,[J]^\perp}$ is essentially unitary, hence Fredholm.

    Let $J=\bigcap\limits_{k=0}^nJ_k$ be an irredundant primary decomposition of $J$ in the sense of Remark \ref{rem:primary_decomposition}, where $J_0$ is finite codimensional in $\poly$, and $\dim_{\mathbb{C}}Z(J_k)\geq1$ for $k\geq1$. Write $s=\mathrm{codim}_{\mathbb{C}}J_0$, then for each positive integer $m$ (to be determined later), if $f$ is a $\mathbf{K}$-quasi-homogeneous polynomial in $\ker S_{z_1^m,[J]^\perp}$ then $z_1^mf\in J\subset J_k,~\forall k\geq1$. By $z_1\notin\sqrt{J_k}$, the definition of primary ideal gives $f\in J_k$, and consequently $f\in\bigcap_{k=1}^nJ_k$. For convenience, write $
L=\bigcap\limits_{k=1}^{n}J_k,
$  and we have $$\ker S_{z_1^m,[J]^\perp}\subset[L]\cap[J]^\perp=[L]\ominus [J].$$
 Since \(J_0\) has finite codimension, $Z(J_0)$ is a finite subset of \(\mathbb D^d\). Moreover, for every
\(\lambda\notin Z(J_0)\), the localizations of \(J\) and \(L\) at
\(\lambda\) coincide. Hence the relative zero set of \(J\) with
respect to \(L\), in the sense of \cite[Theorem 3.1]{Guo1},
is contained in \(Z(J_0)\subset\mathbb D^d\). It follows from  \cite[Theorem 3.1]{Guo1} that the canonical
map
$
\tau: L/J\longrightarrow [L]/[J]
$
is an isomorphism. Therefore,
$$
\begin{aligned}
\dim_{\mathbb C}\ker S_{z_1,[J]^\perp}^{\,m}
&\le \dim_{\mathbb C}\bigl([L]\ominus[J]\bigr)\\
&=\dim_{\mathbb C}L/J\\
&\le\dim_{\mathbb C}\mathbb C[z_1,\cdots,z_d]/J_0:=s.
\end{aligned}
$$
Next, we claim that there is an $m$ such that $\dim_{\mathbb{C}}\ker S_{z_1,[J]^\perp}^{m*}>s$.  To prove the claim, let $\cal{Q}=[J]^\perp$ and consider the $\mathbf{K}$-quasi-homogeneous graded decomposition
$$
H^2(\mathbb D^2)=\bigoplus\limits_{l=0}^\infty {\cal P}_l^{\mathbf{K}},\text{ }[J]=\bigoplus\limits_{l=0}^\infty(J\cap {\cal P}_l^{\mathbf K}),\text{ and } \mathcal{Q}=\bigoplus\limits_{l=0}^\infty\mathcal{Q}_l,$$
where $\mathcal P_{l}^{\mathbf{K}}=\left\{\sum\limits_{\alpha}a_\alpha z^\alpha; \sum\limits_{i=1}^d \alpha_i K_i=l\right\}$ and $\mathcal{Q}_l= {\cal P}_l^{\mathbf{K}}\ominus (J\cap {\cal P}_l^{\mathbf K})
$. Obviously, $S_{z_1}^m \mathcal{Q}_{l}\subset\mathcal{Q}_{l+mK_1}$, and consequently
$$
\mathrm{ran}\, S_{z_1,[J]^\perp}^m\subseteq \bigoplus\limits_{l\geq mK_1}\mathcal{Q}_l.
$$
It follows that
$$\bigoplus\limits_{0\leq l<m K_1}\mathcal{Q}_l\subset (\mathrm{ran}\, S_{z_1,[J]^\perp}^m)^\perp=\ker S_{z_1,[J]^\perp}^{m*}.$$
By the assumption, $$\lim\limits_{r\to\infty}\sum\limits_{l=0}^r \dim\mathcal{Q}_l=\dim\mathcal{Q}=\dim[J]^\perp=\infty.$$ Choose $m$ large enough, such that $\sum\limits_{0\leq l< mK_1}\dim \mathcal Q_l>s$, and hence $\dim\ker S^{m*}_{z_1,[J]^\perp}>s$.
This proves the claim, and moreover $\mathrm{Ind}~S_{z_1,[J]^\perp}^m<0.$ The conclusion of the lemma follows from $$m\mathrm{Ind}~S_{z_1,[J]^\perp}=\mathrm{Ind}~S_{z_1,[J]^\perp}^m.$$
\end{proof}
\begin{thm}\label{thm:K-homology}
    Let $I\subset\poly$ be a Clark-type ideal, then the short exact sequence
    $$0\to \cal{K} \hookrightarrow C^*({[I]^\perp})+\mathcal{K}\to C\big(Z(I)\cap \partial\mathbb D^d\big)\to 0$$
    is not splitting.
\end{thm}
\begin{proof}
    Let $J,B_1,\ldots,B_d$ be as in Definition \ref{def:Clarktype}. By Lemma \ref{lem:nonsplit}, $S_{z_1,[J]^\perp}$ is Fredholm with nonzero index. By \cite[Lemma 5.5]{GWk1}, it suffices to prove that $S_{B_1(z_1),[I]^\perp}$ is Fredholm with nonzero index.

    By the proof of Theorem \ref{thm:enormal_Hardy}, $S_{B_1(z_1),[I]^\perp}$ is unitarily  equivalent to $S_{z_1,[J]^\perp}\otimes\Id_{H_0}$, where the subspace $H_0=\bigotimes_{i=1}^d\ker M_{B_i(z_i)}^*$ is finite dimensional. Since $S_{z_1,[J]^\perp}$ is Fredholm with nonzero index, then so are $S_{z_1,[J]^\perp}\otimes\Id_{H_0}$ and $S_{B_1(z_1),[I]^\perp}$, which completes the proof of the theorem.
\end{proof}

\section{Hilbert-Schmidt submodules vs. Essentially normal quotient modules}

~~~~~~In this section, we establish the connection between Hilbert-Schmidt submodules and essentially normal quotient modules in $\mathcal{A}_\alpha(\mathbb{D}^2)$. Moreover,  Theorems \ref{Theorem:Main} and \ref{Theorem:Main1} will be proved.

Let $p\in\mathbb{C}[z,w]$ be a polynomial such that
$$Z(p)\cap Z(z)\cap\partial\mathbb D^2=\emptyset,$$
i.e. $p(0,w)$ does not vanish in the unit circle, and $\alpha\geq-1$ be an integer. Suppose that $\mathcal{M}$ is a finitely generated submodule of $\mathcal{A}_\alpha(\mathbb{D}^2)$ containing $p$, i.e. $\mathcal{M}=[p,f_1,\cdots, f_n]_\alpha$ where $f_i\in\mathcal{A}_\alpha(\mathbb{D}^2)$, and set
$$Q=\mathcal{M}\ominus[p]_\alpha=[p]_\alpha^\perp\ominus \mathcal{M}^\perp.$$
It is not hard to verify that $Q$ is invariant for $S_{z,[p]_\alpha^\perp}$. In fact, for any $f\in Q\subset\mathcal{M}$, $g\in \mathcal{M}^\perp$ and $q\in\mathbb C[z,w]$, it holds that
$$
\langle S_{q,[p]_\alpha^\perp} f, g\rangle=\langle P_{[p]_\alpha^\perp} M_q f, g\rangle=\langle qf,g\rangle=0.
$$
 For $q\in\mathbb{C}[z,w]$, we write
\begin{eqnarray}\label{eq:TqQ}
T_{q,Q}=S_{q,[p]_\alpha^\perp}\mid_Q=P_QM_q\mid_Q.
\end{eqnarray}
The following lemma is inspired by \cite{LIY}.
\begin{lem}\label{lem:finitecodimension}Under the above assumptions, the subspace $\ker  T_{z,Q}^{*}$ is of finite dimension.
\end{lem}
\begin{proof}
Please note that
\begin{eqnarray*}
\ker T_{z,Q}^{*}&=&Q\ominus\overline{S_{z,[p]_\alpha^\perp} Q}\\
                                                 &=&\{f\in Q:\langle f, P_{[p]_\alpha^\perp}(z g)\rangle=0,~\forall g\in Q\}\\
                                                 &=&\{f\in Q:\langle f, z g\rangle=0,~\forall g\in Q\}\\
                                                 &=&\{f\in Q:\langle f, z g\rangle=0,~\forall g\in\mathcal{M}\}\\
                                                 &=&\{f\in Q:\langle f, h\rangle=0,~\forall h\in z\mathcal{M}\}\\
                                                 &=&\mathcal{M}\cap [p]_\alpha^\perp\cap \left(z\mathcal{M}\right)^\perp\\
                                                 &\subseteq&\mathcal{M}\ominus [z\mathcal{M}+p\cal M]_\alpha.
\end{eqnarray*}

Next, we prove that $\mathcal{M}\ominus[z\mathcal{M}+p\cal {M}]_\alpha$ is of finite dimension. By the assumption
$$
Z(p)\cap Z(z)\cap\partial\mathbb D^2=\emptyset,
$$
we have
$$
\dim_\mathbb{C}\mathbb{C}[z,w]/(z,p)<\infty,
$$
where $(z,p)$ is the ideal generated by $z$ and $p$. Select a basis $\{e_i+(z,p):1\leq i\leq k\}$ of $\mathbb{C}[z,w]/(z,p)$. We claim that
$$\mathcal{M}\ominus [z\mathcal{M}+p\cal M]_\alpha=\mathrm{span}\left\{P(e_i f_j): i=1,\cdots,k,~j=1,\cdots, n+1\right\},$$
where for convenience, we denote $p$ by $f_{n+1}$ and $P$ is the projection onto $\cal M\ominus [z\cal M+p\cal M]_\alpha$. To prove this, suppose $g\in\mathcal{M}\ominus[z\mathcal{M}+p\cal  M]_\alpha$ satisfies
$$g\perp \mathrm{span}\left\{P(e_i f_j): i=1,\cdots,k,~j=1,\cdots, n+1\right\}.$$
Now, for any $h\in\mathbb{C}[z,w]$, there exists $h_1\in\mathrm{span}\{e_1,\ldots,e_k\}$ and $h_2\in(z,p)$, such that $h=h_1+h_2$, then $$P(h_1 f_j)\in\mathrm{span}\{P(e_if_j):~i=1,\ldots,k\}.$$
Furthermore, $h_2 f_j~(1\leq j\leq n+1)$ are in the closure of $z\cal M+p\cal M$.
Therefore
\begin{eqnarray*}
\langle g,h f_j\rangle=\langle g,h_1f_j\rangle+\langle g, h_2 f_j\rangle=\langle g,P(h_1f_j)\rangle+\langle g, h_2 f_j\rangle=0.
\end{eqnarray*}
It follows from the fact $\cal {M}=[p, f_1,\cdots,f_n]$ that $g=0$,  and the claim is proved. Consequently, $\ker  T_{z,Q}^{*}$ is of finite dimension.
\end{proof}


The following lemma originates from Arveson \cite{Arv2,Arv1}, and we omit its proof which is almost trivial.
\begin{lem}\label{lem:decomposition}
Let $A$ be  $(1,\infty)$-essentially normal (respectively, essentially normal) on a Hilbert space, and $\mathcal{M}$ be a closed invariant subspace for $A$. Then the following statements are equivalent.
\begin{itemize}
     \item[(1)] the restriction $A\mid_\mathcal{M}$ is $(1,\infty)$-essentially normal (resp. essentially normal);
     \item[(2)] the operator $P_\mathcal{M}A|_{\mathcal{M}^\perp}$ belongs to $\mathcal{L}^{(2,\infty)}$ (resp. $\cal K$);
     \item[(3)] the compression $P_{\mathcal{M}^\perp}A|_{\mathcal{M}^\perp}$ is $(1,\infty)$-essentially normal (resp. essentially normal).
\end{itemize}
\end{lem}

\begin{thm}\label{thm4.5}
  Let $\mathcal{M}$ be a submodule in $H^2(\mathbb{D}^2)$ containing a Clark-type ideal $I,$ then the
  following statements are equivalent.
\begin{itemize}
     \item[(1)] $\mathcal{M}$ is Hilbert-Schmidt;
     \item[(2)] $\mathcal{M}^\perp$ is $(1,\infty)$-essentially normal;
     \item[(3)] $\mathcal{M}_1=[I]^\perp\ominus\mathcal{M}^\perp$ is $(1,\infty)$-essentially normal.
\end{itemize}
\end{thm}
\begin{proof}
    Theorem \ref{thm:enormal_Hardy} shows that $[I]^\perp$ is $(1,\infty)$-essentially normal. Then the equivalence between (2) and (3) follows immediately from Lemma \ref{lem:decomposition}.

    (1)$\Rightarrow$(2): Suppose $\mathcal{M}$ is Hilbert-Schmidt, i.e.
    $$C_\mathcal{M}=\Id_{\mathcal{M}}-R_zR_z^*-R_wR_w^*+R_zR_wR_z^*R_w^*\in\mathcal{L}^2.$$
    By \cite[Theorem 4.3]{Ya1}, it holds that $\sigma_e(S_z,S_w)\subset\partial\mathbb{D}^2$. On the other hand, for each $f\in I$ we have
    $f(S_z,S_w)=S_f=0$, which forces
    $$f\big(\sigma(S_z,S_w)\big)=\sigma\big(f(S_z,S_w)\big)\subset\{0\}.$$
    Consequently $\sigma_e(S_z,S_w)\subset Z(I)\cap\partial\mathbb{D}^2\subset\mathbb{T}^2$. Applying the Spectral Mapping Theorem, we find $\sigma_e(S_z)\subset\mathbb{T}$ which ensures the Fredholmness of $S_z$. It follows from Theorem \ref{thm:enormal_Hardy} that
    $$\Id_{[I]^\perp}-S_{z,[I]^\perp}S_{z,[I]^\perp}^*\in\mathcal{L}^{(1,\infty)},$$
    and consequently
    $$\Id_{\mathcal{M}^\perp}-S_{z}S_{z}^*\in\mathcal{L}^{(1,\infty)}.$$
    Let $S_z=W|S_z|$ be the polar decomposition, where $|S_z|=(S_z^*S_z)^{1/2}$ and $W$ is the partial isometry that maps $(\ker S_z)^\perp$ unitarily onto $\range S_z$. Then the inclusion
    $$\Id_{\mathcal{M}^\perp}-S_{z}^*S_{z}=W^*(\Id_{\mathcal{M}^\perp}-S_{z}S_{z}^*)W+P_{\ker S_z}\in\mathcal{L}^{(1,\infty)},$$
    derives the $(1,\infty)$-essential normality of $\mathcal{M}^\perp$. 


    (2)$\Rightarrow$(1): This follows from the direct computation
    \begin{eqnarray*}
        C_\mathcal{M}&=&P_{\mathcal{M}}(\Id-M_zP_{\mathcal{M}}M_z^*-M_wP_{\mathcal{M}}M_w^*+M_{zw}P_{\mathcal{M}}M_{zw}^*)\mid_{\mathcal{M}}\\
        &=&P_{\mathcal{M}}(\Id-M_zM_z^*-M_wM_w^*+M_{zw}M_{zw}^*)\mid_{\mathcal{M}}\\
        &&+P_{\mathcal{M}}(M_zP_{\mathcal{M}^\perp}M_z^*+M_wP_{\mathcal{M}^\perp}M_w^*-M_{zw}P_{\mathcal{M}^\perp}M_{zw}^*)\mid_{\mathcal{M}}\\
        &=&P_{\mathcal{M}}C_{H^2(\mathbb{D}^2)}\mid_{\mathcal{M}}+A_zA_z^*+A_wA_w^*-A_{zw}A_{zw}^*\\
        &\in&\mathcal{L}^{(1,\infty)},
    \end{eqnarray*}
    where $A_p=P_{\mathcal{M}}M_pP_{\mathcal{M}^\perp}$ for $p\in\mathbb{C}[z,w]$, which belongs to $\mathcal{L}^{(2,\infty)}$ by Lemma \ref{lem:decomposition}.
\end{proof}
In Example \ref{Exam6.1} in the last section, we will give a counterexample to illustrate that, the submodule $\mathcal{M}$ containing a Clark-type ideal $I$ {\color{red}is not necessarily} Hilbert-Schmidt.

Similar to the submodules in the Hardy space, a submodule ${\cal M}$ in the weighted Bergman module $\mathcal{A}_\alpha(\mathbb{D}^2)$ for integers $\alpha\in [-1,\infty)$ is said to be Hilbert-Schmidt if $C_{\cal M}$, defined in (\ref{Def:Core-Operator}),  is Hilbert-Schmidt.
\begin{rem}
    The proof of (2)$\Rightarrow$(1) in Theorem \ref{thm4.5} is also valid to weighted Bergman modules. In fact, direct computation shows that $C_{\mathcal{A}_\alpha(\mathbb{D}^2)}$ is Hilbert-Schmidt for $\alpha>-1$. Then by Lemma \ref{lem:decomposition}, for a finitely generated  submodule $\cal M$ of $\mathcal{A}_\alpha(\mathbb{D}^2)$ containing a Clark-type ideal, such that $\mathcal{M}^\perp$ is $(1,\infty)$-essentially normal, the core operator
    $$
    C_\mathcal{M}=P_{\mathcal{M}}C_{\mathcal{A}_\alpha(\mathbb{D}^2)}\mid_{\mathcal{M}}+A_zA_z^*+A_wA_w^*-A_{zw}A_{zw}^*
    $$ is Hilbert-Schmidt.
      While the proof of (1)$\Rightarrow$(2) needs \cite[Theorem 4.3]{Ya1}, which was only proven for the Hardy module.
\end{rem}
The following result is the key to prove Theorems \ref{Theorem:Main} and \ref{Theorem:Main1}.

\begin{thm}\label{thm:HSandEM}
    Let $p\in\mathbb{C}[z,w]$ be a distinguished polynomial, and $\alpha\geq-1$ be an integer, such that $S_{z,[p]_\alpha^\perp}$ and $S_{w,[p]_\alpha^\perp}$ are $(1,\infty)$-essentially unitary. If $\mathcal{M}_\alpha$ is a finitely generated submodule of $\mathcal{A}_\alpha(\mathbb{D}^2)$ containing $p$, then the quotient module $\mathcal{M}_\alpha^\perp$ is $(1,\infty)$-essentially normal. As a consequence, the submodule $\cal M_\alpha$ is Hilbert-Schmidt.
\end{thm}
\begin{proof} For $\mathcal{M}_\alpha=[p,f_1,\cdots, f_n]_\alpha,$ we decompose the quotient module $[p]_\alpha^\perp$ into $[p]_\alpha^\perp=Q\oplus \mathcal{M}_\alpha^\perp,$ then accordingly $S_{z,[p]_\alpha^{\perp}}$ has a matrix representation
\begin{eqnarray}\label{eq:decomposition}
S_{z,[p]_\alpha^{\perp}}=\left(\begin{array}{cc} T_{z,Q} & D_z\\ 0 &S_{z,{\mathcal{M}_\alpha^\perp}}\end{array}\right),
\end{eqnarray}
where $D_z=P_{Q}S_z|_{\mathcal{M}_\alpha^\perp}$. By assumption, $S_{z,[p]_\alpha^\perp}$ is $(1,\infty)$-essentially unitary. Since $Q$ is invariant for $S_{z,[p]_\alpha^\perp}$, we conclude that $T_{z,Q}$ is $(1,\infty)$-essentially isometric. Moreover, $\ker T_{z,Q}^*$ is finite dimensional by Lemma \ref{lem:finitecodimension}, therefore $T_{z,Q}$ is Fredholm, hence $(1,\infty)$-essentially unitary. Then by Lemma \ref{lem:decomposition}, $S_{z,\mathcal{M}_\alpha^\perp}=P_{\mathcal{M}_\alpha^\perp}S_{z,[p]_\alpha^{\perp}}\mid_{\mathcal{M}_\alpha^\perp}$ is $(1,\infty)$-essentially normal. Similarly, $S_{w,\mathcal{M}_\alpha^\perp}$ is also $(1,\infty)$-essentially unitary, and hence $\mathcal{M}_\alpha^\perp$ is $(1,\infty)$-essentially normal.
\end{proof}
Obviously,  Theorem \ref{Theorem:Main} is a corollary of Theorem \ref{them:GWZ} and  Theorem \ref{thm:HSandEM}. Moreover, as mentioned in the introduction,  $\mathcal{A}_{-1}(\mathbb D^d)$ is the Hardy space $H^2(\mathbb D^d)$, and hence Theorem \ref{Theorem:Main1} is a corollary of Theorem \ref{thm:enormal} and Theorem \ref{thm:HSandEM}.

\section{Hilbert-Schmidtness of the submodule containing a nonzero polynomial}

~~~~In this section, we study the Hilbert-Schmidtness of finitely generated submodules containing a nonzero polynomial. The following result generalizes Yang's result \cite{Ya3,Ya2} to non-algebraic cases.


\begin{thm}
Let $p$ be a nonzero polynomial, and $\mathcal{M}$ be a finitely generated submodule containing $p$, then  $\mathcal{M}$ is Hilbert-Schmidt.
\end{thm}

\begin{proof} The theorem will be proved in two steps.

\ms
{\bf Step 1.} We first consider the case $Z(p)\cap Z(z)\cap\partial\mathbb D^2=\emptyset$.
\ms
With respect to the decomposition $[p]^\perp=Q\oplus \mathcal{M}^\perp$, $S_{z,{[p]^\perp}}$ can be decomposed as
$$
S_{z,{[p]^\perp}}=\left(\begin{array}{cc}T_{z,Q} & D_z \\ O & S_{z,\mathcal{M}^\perp}\end{array}\right).
$$
We claim that $T_{z,Q}$ is Fredholm. By \cite{Ya2}, the submodule $[p]$ is Hilbert-Schmidt. Then \cite[Theorem 4.3]{Ya1} implies $\sigma_e(S_{z,[p]^\perp},S_{w,[p]^\perp})\subset\partial\mathbb{D}^2$. On the other hand, we have
    $$\sigma_e(S_{z,[p]^\perp},S_{w,[p]^\perp})\subset\sigma(S_{z,[p]^\perp},S_{w,[p]^\perp})\subset Z(p),$$
    and it follows that $\sigma_e(S_{z,[p]^\perp},S_{w,[p]^\perp})\subset Z(p)\cap\partial\mathbb{D}^2$. Since $Z(p)\cap Z(z)\cap\partial\mathbb D^2=\emptyset$, the Spectral Mapping Theorem gives $0\notin\sigma_e(S_{z,{[p]^\perp}})$, hence $S_{z,{[p]^\perp}}$ is Fredholm. Then there is a bounded linear operator $A\in B([p]^\perp$) such that $AS_{z,{[p]^\perp}}\in\Id_{[p]^\perp}+\mathcal{K}$, and consequently $P_{Q}AS_{z,{[p]^\perp}}P_{Q}\in\Id_{Q}+\mathcal{K}$, i.e. $P_{Q}AP_QT_{z,Q}\in\Id_{Q}+\mathcal{K}$. This is equivalent to saying that, $\dim_{\mathbb{C}}\ker T_{z,Q}<\infty$ and $\range T_{z,Q}$ is closed. These facts together with Lemma  \ref{lem:finitecodimension} show the Fredholmness of $T_{z,Q}$.

    Hence there exists a compact operator $K$ and a Fredholm operator $A_z:P_{Q}AP_Q\in B(Q)$ such that $A_zT_{z,Q}=\Id_{Q}+K.$
Obviously,
$
\left(\begin{array}{cc}A_z & O \\ O& \Id_{\mathcal{M}^\perp}\end{array}\right)
$ is Fredholm and so is
$$
\left(\begin{array}{cc}A_z & O \\O & \Id_{\mathcal{M}^\perp}\end{array}\right) \left(\begin{array}{cc}T_{z,Q} & D_z \\O & S_{z,\mathcal{M}^\perp}\end{array}\right)= \left(\begin{array}{cc}\Id_Q+K & D \\O & S_{z,\mathcal{M}^\perp}\end{array}\right),
$$
where $D=A_zD_z$ is a bounded operator. It follows that
$
\left(\begin{array}{cc}\Id_Q & D \\O & S_{z,\mathcal{M}^\perp}\end{array}\right)
$ is Fredholm. Now, from the fact
$$
\left(\begin{array}{cc}\Id_Q & D \\O & \Id_{\mathcal{M}^\perp}\end{array}\right)\left(\begin{array}{cc}\Id_Q & -D \\ O& \Id_{\mathcal{M}^\perp}\end{array}\right)=\left(\begin{array}{cc}\Id_Q & -D \\O & \Id_{\mathcal{M}^\perp}\end{array}\right)\left(\begin{array}{cc}\Id_Q & D \\O & \Id_{\mathcal{M}^\perp}\end{array}\right)=\left(\begin{array}{cc}\Id_Q & O\\O & \Id_{\mathcal{M}^\perp}\end{array}\right),
$$
we conclude that $\left(\begin{array}{cc}\Id_Q & -D \\O & \Id_{\mathcal{M}^\perp}\end{array}\right)$ is invertible, and hence
$$
\left(\begin{array}{cc}\Id_Q & O \\O & S_{z,\mathcal{M}^\perp}\end{array}\right)=\left(\begin{array}{cc}\Id_Q & D \\O & S_{z,\mathcal{M}^\perp}\end{array}\right)\left(\begin{array}{cc}\Id_Q & -D \\O & \Id_{\mathcal{M}^\perp}\end{array}\right)
$$ is Fredholm, which ensures the Fredholmness of $S_{z,\mathcal{M}^\perp}$, i.e. $0\not\in\sigma_e(S_{z,\mathcal{M}^\perp})$. Then by \cite[Theorem 2.3]{Ya2}, $\mathcal{M}$ is Hilbert-Schmidt.
\ms

{\bf Step 2.} For the case $Z(p)\cap Z(z)\cap\partial\mathbb D^2\neq\emptyset$, write
$$
p=\prod\limits_{j=1}^n (w-\xi_j)^{k_j} p_0,
$$
where $|\xi_j|=1$. Obviously, $[p]=[p_0]$. Without loss of generality, we may assume that $p=p_0$. Now, if $p_0$ is a nonzero constant, then $\cal {M}=[p]=H^2(\mathbb D^2)$ and obviously $\cal M$ is Hilbert-Schmidt. Next, assume $p$ is a polynomial such that $Z(p)\cap \mathbb D^2\not=\emptyset.$ By \cite[Corollary 3.5]{Ya3}, there is an $a\in\mathbb D$ such that $p_0$ has no zeros in $\{a\}\times\mathbb T$. Denote the M\"obius transformation by
$$
\varphi_a(z)=\frac{a-z}{1-\bar{a}z}, \quad \Phi_a(z,w)=\big(\varphi_{a}(z), w\big).
$$
It is easy to see that
$$
Z(z)\cap Z(p_{0}\circ\Phi_a)\cap\partial \mathbb D^2=\emptyset.
$$
Let $\mathcal{M}_{\Phi_a}=\{f\circ \Phi_a: f\in \mathcal M\}$,
then by Step 1  the core operator $C_{\mathcal{M}_{\Phi_a}}$ is Hilbert-Schmidt.  Then the Hilbert-Schmidtness of $C_\mathcal{M}$ follows from \cite[Proposition 2.1]{Ya1}.
\end{proof}

\section{Examples}

~~~~In this section, we consider several examples. At first, we show that in Theorem \ref{thm:HSandEM}, without the assumption that the submodule is finitely generated, the associated quotient module may fail to be essentially normal.
\begin{exam}\label{Exam6.1}
Consider the quotient module $[z-w]^\perp$ of $H^2(\mathbb D^2)$, which can be seen as the Bergman space over the unit disc. Take an invariant subspace $\mathcal{N}$ of $S_z$ such that $\dim(\mathcal{N}\ominus S_z\mathcal{N})=\infty$, then the submodule $\mathcal{M}=[z-w]\oplus\mathcal{N}$ is infinitely generated. By \cite{LIY}, $\mathcal{M}$ is not Hilbert-Schmidt, and hence by Theorem \ref{thm4.5}, $\mathcal{M}^\perp$ is not essentially normal.
\end{exam}
\hfill $\square$

We conjecture that all the quotient modules corresponding to infinitely generated submodules are not
essentially normal. In the following example, we prove that this is the case for Rudin's quotient module.
\begin{exam}
In \cite[Page71]{Ru}, Rudin constructed an invariant subspace of $H^2(\mathbb{D}^2)$ of infinite rank as follows. It is the submodule $\mathcal{M}$ consisting of all functions in $H^2(\mathbb{D}^2)$ with a zero of order at least $n$ at $\left(\alpha_n, 0\right)=\left(1-n^{-3}, 0\right)$ for each positive integer $n$. Let
$$
\left\{\begin{array}{l}
b_n(z)=\left(z-\alpha_n\right) /\left(1-\alpha_n z\right),\\
q_0(z)=\prod_{n=1}^{\infty} b_n^n(z), \\
q_j(z)=q_{j-1}(z) / \prod_{n=j}^{\infty} b_n(z),~j \geq 1.
\end{array}\right.
$$
In \cite{SY}, Seto and Yang obtained the decompositon
$$\mathcal{M}^\bot=\bigoplus_{j=0}^{\infty}\big( H^2(z) \ominus q_j(z) H^2(z)\big) w^j.$$
In the case $ n\geq 2,$ since $\alpha_n\in Z(q_0)\cap Z(q_1),$
it is easy to see that the normalized reproducing kernel
$$k_{\alpha_n}(z) \in H^2(z) \ominus q_1(z) H^2(z)\subseteq H^2(z) \ominus q_0(z) H^2(z).$$
Routine calculation gives
$$(S_w^*S_w-S_wS_w^*)k_{\alpha_n}(z)=k_{\alpha_n}(z),$$
therefore $S_w^*S_w-S_wS_w^* $ is not compact, and consequently $\mathcal{M}^\bot$ is not essentially normal.
\end{exam}
\hfill $\square$

As usual, for a quotient module ${\mathcal M}^\perp$ in $H^2(\mathbb D^2)$, we write $\sigma_e(\mathcal{M}^\perp)$ for $\sigma_e(S_{z,{\mathcal M}^\perp},S_{w,{\mathcal M}^\perp})$.

\begin{exam}
In this example, we construct a non-algebraic submodule $\mathcal{M}$ containing $p_B$ such that $\mathcal{M}\ominus[p_B]$ is infinite dimensional, while $\sigma_e(\mathcal{M}^{\perp})=\sigma_e([p_B]^{\perp}).$

Let $p=\prod\limits_{i=1}^n\left(z-\alpha_i w\right)^{n_i}$ be a distinguished homogeneous polynomial where each $|\alpha_i|=1,$ and $B_1,B_2$ be two nonconstant finite Blaschke products. As before, set $B=(B_1, B_2)$ and
$$
p_{B}(z,w)=\prod\limits_{i=1}^n\left(B_1(z)-\alpha_i B_2(w)\right)^{n_i}.
$$  Let $\tilde{B}$ be an infinite Blaschke product satisfying $\mathbb{T} \subseteq \overline{Z(\tilde{B})}.$ Set $D=Z(\tilde{B}\circ B_1)$, and
$$E=\bigcup_{i=1}^n Z\big(\tilde{B}\circ(\alpha_iB_2)\big)=Z\Big(\prod_{i=1}^n\tilde{B}\circ(\alpha_iB_2)\Big).$$
Set $B_3=\tilde{B}\circ B_1$ and $B_4=\prod_{i=1}^n\tilde{B}\circ(\alpha_iB_2)$, then by Theorem \ref{thm:HSandEM} and \cite[Theorem 4.3]{Ya1}, the submodule $\mathcal{M}=[B_3(z)B_4(w),p_B]$ is Hilbert-Schmidt, and moreover, $\sigma_e(\mathcal{M}^{\perp}) \subseteq \partial \mathbb{D}^2$.
Then the Spectral Mapping Theorem ensures that
$$\sigma_e(\mathcal{M}^{\perp}) \subseteq Z(p_B) \cap \partial \mathbb{D}^2=\sigma_e([p_B]^\perp).$$

Conversely, for $(\lambda_1,\lambda_2)\in\sigma_e([p_B]^\perp),$ since
$B_1(\lambda_1)=\alpha_iB_2(\lambda_2)\in \mathbb{T}$ for some $i,$ there exists a sequence $\{\mu_n\}$ in $Z(\tilde{B})$ such that $\mu_n\to \tau$, where $$\tau=B_1(\lambda_1)=\alpha_iB_2(\lambda_2)\in\mathbb T.$$ Then by \cite[Lemma 2.2]{DPW}, we can find $t_n\in\{B_1^{-1}(\mu_n)\}\subset Z(B_3)$ and $t_n^{\prime}\in \{(\alpha_iB_2)^{-1}(\mu_n)\}\subset Z(B_4)$
such that $(t_n, t_n^{\prime})\rightarrow (\lambda_1,\lambda_2).$ Observing that
$$p_B(t_n,t_n')=p(\eta_n,\alpha_i^{-1}\eta_n)=0,$$
we conclude
$$(t_n,t_n')\in Z\big(B_3(z)B_4(w)\big)\cap Z(p_B).$$
Then by \cite[Lemma 2.2]{DPW} and \cite[Theorem 6.1]{GWp1}, it is not difficult to verify
$$\sigma_e([p_B]^\perp)\subseteq\overline{Z\big(B_3(z)B_4(w)\big)\cap Z(p_B)}\cap\partial\mathbb{D}^2\subseteq\sigma_e(\mathcal{M}^{\perp}).$$

Finally, notice that if $\cal M$ were algebraic, then $\dim Z(\cal M)=1$ or $\dim Z(\cal M)$ is a finite set. However $Z(\cal M)$ is an infinite, discrete set, and it follows that $\cal M$ is not algebraic.

\end{exam}
\hfill $\square$

 By Theorem \ref{thm:HSandEM}, both $\mathcal{M}^\bot=\left[B_3(z)B_4(w),p_B\right]^\bot$ and $\left[ p_B\right]^\bot$ are essentially normal,
hence there are two $C^*$-extensions
\begin{eqnarray}\label{Short-Exact-1}
0\rightarrow \mathcal{K} \rightarrow C^*(M^\bot)+\cal {K}\rightarrow C\big(\sigma_e([p_B]^\bot)\big) \rightarrow  0,\end{eqnarray}
and
\begin{eqnarray}\label{eq:exactsequence} 0\rightarrow \mathcal{K} \rightarrow C^*([p_B]^\bot)+\mathcal K\rightarrow C\big(\sigma_e([p_B]^\bot)\big) \rightarrow  0.
\end{eqnarray}
From Theorem \ref{thm:K-homology}, the short exact sequence (\ref{eq:exactsequence}) gives a nontrivial element in the K-homology group $K\big(\sigma_e([p_B]^\bot)\big)$.
Naturally, one has the following question.
\begin{ques}
Are these two $C^*$-extensions essentially unitarily equivalent?
\end{ques}
The answer to this question is negative, in fact we have the following Theorem.
\begin{thm}\label{Them:Ext}
Let $[e_{[p_B]^\perp}]$ and $[e_{M^\perp}]$ be the K-homology element in $\Ext\big(C\big(\sigma_{e}([p_B]^\perp)\big)\big)$, defined by the short exact sequences  (\ref{Short-Exact-1}) and (\ref{eq:exactsequence}) separately, where $\Ext\big(C\big(\sigma_{e}([p_B]^\perp)\big)\big)$ is the extension group. Then $[e_{[p_B]^\perp}]\not=[e_{M^\perp}]$.
\end{thm}
To prove this theorem, we need the following key technical lemma.
\begin{lem}\label{lem:invert}
There is a $\lambda\in\mathbb D$, such that $S_{B_1, {\cal M}^\perp}-\lambda I$ is invertible.
\end{lem}
\begin{proof}
Since $Z(\tilde{B})$ is countable, we can take $\lambda\in \mathbb D$ such that
\begin{eqnarray}\label{eq:ZeroB}
\tilde{B}(\alpha_i\alpha^{-1}\lambda)\not=0, \text{ for all } 1\leq i,j\leq n.
\end{eqnarray}
Define $F(u,v)={\color{red}\tilde{B}}\prod\limits_{j=1}^n\tilde{B}(\alpha_j v)$ and set $g(z,w)=B_3(z)B_4(w)$. Then
$$
g(z,w)=F\big(B_1(z),B_2(w)\big).
$$
Consider the polynomial $p_\lambda$ on $v$ defined by
$$
p_\lambda(v)=\prod\limits_{i=1}^n(\lambda-\alpha_i v)^{n_i},
$$
then $z(p_\lambda)=\{\alpha_i^{-1}\lambda\}_{i=1}^n$, and $\ord_{p_\lambda}(\alpha_i^{-1}\lambda)=n_i$, where $\ord_{f}(\lambda_0)$ is the order of the zero $\lambda_0$ of $f$. By the choice of $\lambda$ (\ref{eq:ZeroB})
$$
F(\lambda,\alpha_i^{-1}\lambda)=\tilde{B}(\lambda)\prod\limits_{j=1}^n(\alpha_j\alpha_i^{-1}\lambda)\not=0,
$$
then $Z(p_\lambda)\cap Z(F_\lambda)=\emptyset,$ where $F_\lambda(\cdot)=F(\lambda,\cdot)$. By the Hermite  interpolation by polynomials on one variable, we have a polynomial $c(\cdot)$ on $v$ such that
$$
p_\lambda(\cdot)\big|{\color{red}(1-c(\cdot)F_\lambda(\cdot)),}
$$
and
\begin{eqnarray}\label{OrderZero}
\ord_{p_\lambda}(v_i)\leq \ord_{1-c(\cdot)F_\lambda(\cdot)}(v_i).
\end{eqnarray}
Therefore $a(v):=\frac{1-c(v)F(\lambda,v)}{p(\lambda, v)}\in H^\infty(\mathbb D)$, and we have
$$
1=c(\cdot)F_\lambda(\cdot)+a(\cdot) p_\lambda(\cdot).
$$
Now, set
$$
d(u,v)=\frac{1-c(v)F(u,v)-a(v)p(u,v)}{u-\lambda}.
$$
By (\ref{OrderZero}), $d$ is holomorphic. Furthermore,
\begin{eqnarray*}
  d(u,v) &=& -c(v)\frac{F(u,v)-F(\lambda,v)}{u-\lambda}-a(v)\frac{p(u,v)-p(\lambda,v)}{u-\lambda} \\
         &=&-c(v)\frac{{\color{red}\tilde(B)}(u)-\tilde{B}(\lambda)}{u-\lambda}\prod\limits_{j=1}^n \tilde{B}(\alpha_j v)-a(v)\frac{p(u,v)-p(\lambda,v)}{u-\lambda}.
\end{eqnarray*}
It follows that $d(\cdot,\cdot)\in H^\infty(\mathbb D^2)$. Obviously,
\begin{eqnarray}\label{equation1}
1=c(v)F(u,v)+a(v)p(u,v)+(u-\lambda)d(u,v).
\end{eqnarray}
Substituting $u=B_1(z)$ and $v=B_2(w)$ in (\ref{equation1}), we have
\begin{eqnarray}
  c\big(B_2(w)\big)g(z,w)+a\big(B_2(w)\big) p_{B}(z,w)+\big(B_1(z)-\lambda\big)d\big(B_1(z),B_2(w)\big)=1.
\end{eqnarray}
Now, notice that $c\big(B_2(w)\big)g(z,w)+a\big(B_2(w)\big) p_{B}(z,w)\in \cal M$, we have that $$S_{c(B_2(w))g(z,w)+a(B_2(w)) p_{B}(z,w)}=0,$$ and hence
$$
(S_{B_1(z)}-\lambda I) S_{d(B_1(z),B_2(w))}=S_{d(B_1(z),B_2(w))}(S_{B_1(z)}-\lambda I)=I.
$$
\end{proof}

Now, we are going to prove Theorem \ref{Them:Ext}.

\vskip1.5mm\noindent{\emph{ Proof of Theorem \ref{Them:Ext}}.} By the proof of Theorem \ref{thm:K-homology},
$$
\ind(S_{B_1, [p_B]^\perp})\not=0.
$$
Moreover, Since $S_{B_1, [p_B]^\perp}$ is an essentially unitary  operator, we have that $\sigma_{e}(S_{B_1, [p_B]^\perp})\subseteq\mathbb T$. It follows that $\{S_{B_1, [p_B]^\perp}-t\lambda I; 0\leq t\leq 1\}$ is a family of Fredholm operators. It follows that
$$
\ind(S_{B_1, [p_B]^\perp}-\lambda I)=\ind(S_{B_1, [p_B]^\perp})\not=0.
$$
On the other hand, by Lemma \ref{lem:invert},
$$
\ind(S_{B_1, {\cal M}^\perp}-\lambda I)=0.
$$
By BDF theory\cite{BDF},
$$
[e_{[p_B]^\perp}]\not=[e_{M^\perp}].
$$

\section{Appendix}
In this section, let $\theta, \psi$ be two nonconstant inner functions, ${\cal M}_{N; \theta, \psi}=\big[\big(\theta(z)-\psi(w)\big)^N\big]$. We derive a trace formula for the  Schatten-$p$ norm of the core operator ${C}_{\mathcal M_{N;\theta,\psi}}$. The key idea is to calculate the Schatten-$p$ norm of  the core operator of the submodule $\big[\big(z-\varphi_{r}(w)\big)^N\big]$, where $\varphi_a=\frac{a-\xi}{1-{\bar a}\xi}$ for $a, \xi\in \mathbb D$ is the M$\ddot{\text{o}}$bius transformation on $\mathbb D$.

For a submodule $\cal M\subset H^2(\mathbb D^2)$, the core operator $C_{\cal M}$ is extended to  an operator on $H^2(\mathbb D^2)$ by
$$
\mathfrak{C}_{M}=\left(\begin{array}{cc} C_{\cal M} & \\ & 0\end{array}\right): {\cal M}\oplus {\cal M}^\perp \to {\cal M}\oplus {\cal M}^\perp.
$$
Obviously, $C_{\cal M}$ and $\mathfrak{C}_{\cal M}$ have the same nonzero spectrum and the same nonzero singular values, and hence they have the same Schatten-$p$ norm.

Next, for inner functions $\theta,\psi$, define the composition operator
$$
({\cal C}_{\theta,\psi}f)(z,w)=f\big(\theta(z),\psi(w)\big),
$$
and let ${\cal M}_{\theta,\psi}$ be the submodule generated by ${\cal C}_{\theta,\psi}\cal M$, i.e.
$$
{\cal M}_{\theta,\psi}=[{\cal C}_{\theta,\psi}\cal M].
$$
It should be pointed out that, in general, ${\cal C}_{\theta,\psi}\cal M$ is not a submodule in $H^2(\mathbb D^2)$. The following lemmas are the key observations.
\begin{lem}\label{lem:Composition}
For any nonconstant inner functions $\theta,\psi$ and submodule $\cal M\subset H^2(\mathbb D^2)$,
$$
\mathfrak{C}_{{\cal M}_{\theta,\psi}}={\cal C}_{\theta,\psi} \mathfrak{C}_{\cal M}{\cal C}_{\theta,\psi}^*.
$$
\end{lem}
\begin{proof}
Let $K_{\lambda,\mu}(z,w)=\frac{1}{(1-\bar\lambda z)(1-\bar\mu w)}$ be the reproducing kernel of $H^2(\mathbb D^2)$. Since for any $(\lambda_i,\mu_i)\in\mathbb D^2$, $i=1,2$,
\begin{eqnarray*}
\langle\cal{C}_{\theta,\psi}^*K_{\lambda_1,\mu_1},K_{\lambda_2,\mu_2}\rangle= \langle K_{\lambda_1,\mu_1},\cal{C}_{\theta,\psi}K_{\lambda_2,\mu_2}\rangle
=\overline{K_{\lambda_2,\mu_2}\big(\theta(\lambda_1),\psi(\mu_1)\big)}=\langle K_{\theta(\lambda_1),\psi(\mu_1)}, K_{\lambda_2,\mu_2}\rangle,
\end{eqnarray*}
we have
$$
\cal{C}_{\theta,\psi}^*K_{\lambda,\mu}=K_{\theta(\lambda),\psi(\mu)}.
$$
Let $K_{\lambda,\mu}^{\cal M}=P_{\cal M}K_{\lambda,\mu}$ be the reproducing kernel of $\cal M$. By \cite{GY}, the core function is defined by
$$
G_{\lambda,\mu}^{\cal M}:=\mathfrak{C}_{\cal M}{K_{\lambda,\mu}}=(1-\bar\lambda z)(1-\bar\mu w)K_{\lambda,\mu}^{\cal M}(z,w).
$$
By \cite[Theorem 3.1]{ZL},
$$
K_{\lambda,\mu}^{M_{\theta,\psi}}(z,w)=\frac{1-\overline{\theta(\lambda)}\theta(z)}{1-\bar\lambda z}\frac{1-\overline{\psi(\mu)}\psi(w)}{1-\bar\mu w} K_{\theta(\lambda),\psi(\mu)}^{\mathcal M}\big(\theta(z),\psi(w)\big),
$$
which is equivalent to
$$
G_{\lambda,\mu}^{{\cal M}_{\theta,\psi}}(z,w)=G_{\theta(\lambda),\psi(\mu)}^{\cal M}\big(\theta(z),\psi(w)\big).
$$
It follows that for any $(\lambda, \mu)\in\mathbb D^2$
\begin{eqnarray*}
{\cal C}_{\theta,\psi} \mathfrak{C}_{\cal M}{\cal C}_{\theta,\psi}^* K_{\lambda,\mu}&=&{\cal C}_{\theta,\psi} G_{\theta(\lambda),\psi(\mu)}^{\cal M}\\
&=&G_{\theta(\lambda),\psi(\mu)}^{\cal M}\big(\theta(z),\psi(w)\big)\\
&=& G_{\lambda,\mu}^{{\cal M}_{\theta,\psi}}(z,w)\\
&=& \mathfrak{C}_{{\cal M}_{\theta,\psi}}K_{\lambda,\mu},
\end{eqnarray*}
and we have the desired result.
\end{proof}
The following easy result is standard, and we write it here for convenience.
\begin{lem}\label{lem:isometry}
Let $\beta$ be an inner function such that $\beta(0)=0$, then the composition operator is an isometry on $H^2(\mathbb D)$.
\end{lem}
\begin{proof}
It is easy to see that $\{\beta^m\}_{m=0}^\infty$ is an orthonormal set in $H^2(\mathbb D)$. In fact, for $m>n\geq 0$,
$$
\langle \beta^m,\beta^n\rangle=\int_{\mathbb T}\beta^{m-n}dm=\beta(0)^{m-n}=0.
$$
Now, for any $f(z)=\sum\limits_{n=0}^\infty a_n z^n\in H^2(\mathbb D)$, we have
$$
\|{\cal C}_\beta f\|^2=\left\|\sum\limits_{n=0}^\infty a_n \beta^n\right\|^2=\|f\|^2.
$$
\end{proof}
As an easy application, for nonconstant inner functions $\theta,\psi$ such that $\theta(0)=\psi(0)=0$,  the composition operator
$$
\mathcal{C}_{\theta,\psi} =\mathcal{C}_{\theta}\otimes\mathcal{C}_{\psi},
$$
and hence $\mathcal{C}_{\theta,\psi}$ is an isometry.

Next, for $a\in \mathbb D$, notice that
$$
\big(\varphi_a(\theta)-\varphi_a(\psi)\big)^N=\left(\frac{1-|a|^2(\psi-\theta)}{(1-\bar{a}\theta)(1-\bar a \psi)}\right)^N.
$$
Since both $1-\bar{a}\theta$ and $1-\bar{a}\psi$ are invertible elements in $H^\infty(\mathbb D^2)$, we have that
\begin{eqnarray}\label{eq:moduletransfer}
\big[\big(\varphi_a(\theta)-\varphi_a(\psi)\big)^N\big]=\big[\big(\theta-\psi\big)^N\big].
\end{eqnarray}
Next, set
$$
b_1=\theta(0),\,b_2=\psi(0), \text{ and } r=\rho_{\mathbb D}(b_1,b_2)=|\varphi_{b_1}(b_2)|,
$$
where $\rho_{\mathbb D}(b_1,b_2)$ is the pseudohyperbolic distance between $b_1$ and $b_2$.

Taking $a=b_1$ in (\ref{eq:moduletransfer}) and then multiplying by a unimodular constant, we may assume without loss of generality that
$$
\theta(0)=0 \text{ and } \psi(0)=r\in[0,1).
$$
For simplicity, set $\varrho_r=\varphi_r\circ \psi$, then $\varrho_r$ is an inner function such that $\varrho_r(0)=0$. Moreover, $\psi=\varphi_r\circ \varrho_r$, we have
$$
\mathcal{C}_{\theta,\psi}=\mathcal{C}_{id, \varrho_r}\mathcal{C}_{\theta,\varphi_r}.
$$
By Lemma \ref{lem:Composition} and Lemma \ref{lem:isometry}, the isometry $V=\mathcal{C}_{\theta,\varrho_r}$ satisfies
$$
\mathfrak{C}_{\mathcal{M}_{N;\theta,\psi}}=V \mathfrak{C}_{\mathcal{M}_{N;z,\varphi_r(w)}}V^*,
$$
where $\mathcal{M}_{N;z,\varphi_r(w)}=\big[\big(z-\varphi_r(w)\big)^N\big]$. We summarize this argument as the following proposition.
\begin{prop}\label{Prop:trans}
For nonconstant inner functions $\theta$ and $\psi$, let $r=\big|\varphi_{\theta(0)}\big(\psi(0)\big)\big|$, then $C_{\mathcal{M}_{N;\theta,\psi}}$ is in the Schatten-$p$ class $\mathcal{S}_p$ for some $p\geq 1$  if and only if $C_{\mathcal M_{N;\varphi_r}}$ is in $\mathcal{S}_p$; moreover, in this case
$$
\|\mathcal{C}_{N;\theta,\psi}\|_p=\|\mathcal{C}_{N;z,\varphi_r}\|_p.
$$
\end{prop}
In the following, we will write $$\mathcal M_{N;\varphi_r}=\big[\big(z-\varphi_r(w)\big)^N\big] \text{ and } \mathcal M_N=\big[\big(z-w)^N].$$ Very recently, Lu and Zu\cite[Theorem 5.1]{LZ} proved that
$$
\|C_{\mathcal M_N}\|_p^p=1+2N^p\zeta(p, N+1),
$$
where Hurwitz zeta function is defined by  $$\zeta(p,N+1)=\sum\limits_{n=1}^\infty\left(\frac{1}{N+n}\right)^p.$$
Next, we are going to deal with the case $r\not=0$. To do this, consider the unitary operator on $H^2(\mathbb D)$
$$
U_r f=f\circ \varphi_r\cdot k_r
$$
where $k_r(z)=\frac{\sqrt{1-r^2}}{1-r z}$ is the normalized reproducing kernel of $H^2(\mathbb D)$ at $r\in\mathbb [0,1)$ and let
$$\mathbb U_r=I\otimes U_r.$$  Set $R_{z}=M_z|_{\mathcal M_N},\, R_{w}=M_w|_{\mathcal M_N}$ temporarily. Denote by
$$
V_{r,w}=\varphi_r(R_w)=(rI-R_w)(I-rR_w)^{-1}.
$$
It can be seen that $V_{r,w}$ is an isometry.
By direct calculations,
$$
\mathbb U_r^* \mathfrak{C}_{\mathcal M_{N;\varphi_r}} \mathbb U_r=I-R_zR_z^*-V_{r,w}V_{r,w}^*+R_zV_{r,w}R_z^*V_{r,w}^*.
$$
Write $P=I-R_zR_z^*$, which is the orthogonal projection onto
\begin{eqnarray}
W_N=\mathcal M_N\ominus z\mathcal M_N
\end{eqnarray}
 then
\begin{eqnarray}\label{eq:6.2}
\mathbb U_r^* \mathfrak{C}_{\mathcal M_{N;\varphi_r}} \mathbb U_r=P-V_{r,w}PV_{r,w}^*.
\end{eqnarray}
Furthermore, $V_{r,w} P V_{r,w}^*$ is the orthogonal projection onto $V_{r,w}W_N$. Now, set
$$
T_N=P R_w|_{W_N}.
$$
Obviously, $R_{w}^*W_N\subset W_N$. For $n>0$, $T_N^n=PR_w^n|_{W_N}$. Moreover,
$$
PV_{r,w}|_{W_N}=P\varphi_r(R_w)|_{W_N}=\varphi_r(T_N).
$$

For convenience, let $\mathbb N_0$ be the nonnegative integers. The following weighted-shift model is contained in \cite{LZ} and we include the argument for completeness.

\begin{prop}
\label{prop:jacobi-model}
There is an orthonormal basis $\{e_k\}_{k\in\mathbb N_0}$ of $W_N$ such that $T_N$ can be represented as the following weighted shift:
$$
 T_Ne_k=\alpha_{N,k}e_{k+1},
$$
where the weight
$$
\alpha_{N,k}=\frac{\sqrt{(k+1)(k+2N+1)}}{k+N+1}.
$$
\end{prop}

\begin{proof}
For $k\in\mathbb N_0$, consider the $N+k$-order homogeneous components of $\mathcal M_N$ and $z{\mathcal M_N}$ separately
$$
\mathcal{M}_{N, N+k}=\big\{(z-w)^Nz^{\alpha_1} w^{\alpha_2};\, (\alpha_1,\alpha_2)\in\mathbb N_0\times\mathbb N_0 \text{ and } \alpha_1+\alpha_2=k\big\},
$$
and
$$
\left(z\mathcal{M}_N\right)_{N+k}=z\mathcal M_{N, N+k-1}=\big\{z(z-w)^Nz^{\alpha_1} w^{\alpha_2};\, (\alpha_1,\alpha_2)\in\mathbb N_0\times\mathbb N_0 \text{ and } \alpha_1+\alpha_2=k-1\big\}.
$$
It follows that $\left(z\mathcal{M}_N\right)_{N+k}$ is a 1-codimensional subspace of $\mathcal{M}_{N,N+k}$, i.e.
$$
W_{N,N+k}:=\mathcal{M}_{N,N+k}\ominus \left(z\mathcal{M}_N\right)_{N+k}=\mathbb C\xi_{N,k},
$$
 and hence $T_N e_{k}=\alpha_{N,k} e_{k+1}$ with the weights $\{\alpha_{N, k}\}$.

 Next, we will use the standard  argument from the theory of orthogonal polynomial on the unit circle \cite{SimonOPUC} to determine the weight $\{\alpha_{N,k}\}$. Now, let
  $$
 \mathcal P_{k}=\big\{f\in\mathbb C[\zeta]; \deg(f)\leq k\big\},
 $$
    where $\mathbb C[\zeta]$ is the polynomial ring on $\zeta$. To identify its canonical generator, for $(z,w)\in\mathbb T^2$ write $\zeta=z/w$. Then the map
$
\Psi: w^{N+k}(\zeta-1)^Nq(\zeta)\mapsto q(\zeta)
$
identifies $\mathcal{M}_{N,N+k}$ with $\mathcal P_{k}$, where circular-Jacobi inner product on $\mathcal P_{k}$ is given by
$$
 \langle f,g\rangle_N
 =\int_{\mathbb T}f(\zeta)\overline{g(\zeta)}|1-\zeta|^{2N}\,dm(\zeta).
$$
Under this identification, $(z\mathcal{M}_{N})_{N+k}$ corresponds to
$\zeta\mathcal P_{k-1}$.  Let $\Phi_{N,k}$ be the monic orthogonal polynomial of degree $k$ for the weight $|1-\zeta|^{2N}$, and its reversed polynomial is
$$
 q_{N,k}(\zeta)=\Phi_{N,k}^*(\zeta)
 =\zeta^k\overline{\Phi_{N,k}(\zeta)},\, |\zeta|=1.
$$
By using standard OPUC theory, we have
$
 q_{N,k}\perp \zeta^j,\, j=1,\cdots, k.
$
Therefore
$$
\mathcal{M}_{N,N+k}=\text{span}\left\{\xi_{N,k}(z,w)=(z-w)^N w^k q_{N,k}(z/w);\,  q_{N,k}\in \mathcal{P}_{k}\right\}.
$$
By  \cite[Chapter~1]{SimonOPUC}, for the circular-Jacobi weight $|1-\zeta|^{2N}$, the Verblunsky coefficients are
$$
 \alpha_j=-\frac{N}{N+j+1},\, j\ge0.
$$
By the direct computation,
$$
 \int_{\mathbb T}|1-\zeta|^{2N}\,dm(\zeta)=\binom{2N}{N}.
$$
Therefore,
\begin{align*}
 \|q_{N,k}\|_N^2 =\binom{2N}{N}  \prod_{j=0}^{k-1}\left(1-\frac{N^2}{(N+j+1)^2}\right)=\frac{(2N+k)!\,k!}{(N+k)!^2}=h_{N,k}.
\end{align*}
Since the constant terms of $q_{N,k}$ and $q_{N,k+1}$ are both $1$, so
$q_{N,k+1}-q_{N,k}\in\zeta\mathcal P_k$. It follows from $q_{N,k+1}\perp \zeta\mathcal P_k$ that $$\langle q_{N,k+1},q_{N,k+1}-q_{N,k}\rangle=0,$$
and hence
$$
 \langle q_{N,k},q_{N,k+1}\rangle_N=
\langle q_{N,k}-q_{N,k+1}, q_{N,k+1}\rangle+\langle q_{N,k+1}, q_{N,k+1}\rangle
 =\|q_{N,k+1}\|_N^2=h_{N,k+1}.
$$
Set
$e_k=\xi_{N,k}/\sqrt{h_{N,k}}$, and  we have
$$
 \langle T_Ne_k,e_{k+1}\rangle
 =\frac{h_{N,k+1}}{\sqrt{h_{N,k}h_{N,k+1}}}
 =\sqrt{\frac{h_{N,k+1}}{h_{N,k}}}.
$$
It follows that
$$\alpha_{N,k}=\sqrt{h_{N,k+1}/h_{N,k}}=\frac{\sqrt{(k+1)(k+2N+1)}}{k+N+1}.
$$
\end{proof}

The following lemma is a direct specialization of Halmos two-projections theorem \cite{Halmos}, which is also a special case of the general formula of Omladi\v{c} \cite{Omladic}. We list the proof here for reader's convenience.

\begin{lem}
\label{lem:two-projections}
Let $P_W$ be the orthogonal projection onto a closed subspace $\mathcal W$ of a Hilbert space $\mathcal H$,  $V$ be an isometry. Let
$$
 Q=VPV^*,\, A=P_{\mathcal W}V|_{\mathcal W},\, D=I_{\mathcal W}-AA^*.
$$
Assume that $A$ is injective, $\dim\ker A^*=1$, and $D$ is compact.  Assume that$$\sigma_p(D)\setminus\{0\}={\color{red}\{1,\lambda_k, k=1,\cdots\}}.$$
Rearranged the nonzero eigenvalue of $D$, counting multiplicities, as
$
\lambda_1\leq\lambda_2\leq\cdots<1.
$
Then
\begin{eqnarray}\label{eigenvalue}
 \sigma(P_{\mathcal W}-Q)\setminus\{0\}
 =\{1\}\cup\{\pm\sqrt{\lambda_j}:j\ge1\}.
\end{eqnarray}
Consequently, for every $0<p<\infty$,
$$
 \|P_{{\mathcal W}}-Q\|_{\mathcal S_p}^p
 =2\Tr(D^{p/2})-1,
$$
where sides may be $+\infty$,

\end{lem}

\begin{proof}
According to the orthogonal decomposition $\mathcal H=\mathcal W\oplus\mathcal W^\perp$,  we have a block matrix
$$
 V|_{\mathcal W}=\binom{A}{B}.
$$
Because $V|_{\mathcal W}$ is an isometry, i.e. $V^*_{\mathcal W}V_{\mathcal W}=I$, we have
$$
 A^*A+B^*B=I_{\mathcal W},\text{ and }
 Q=VP_{\mathcal W}V^*
 =\begin{pmatrix}
 AA^*&AB^*\\ BA^*&BB^*
 \end{pmatrix}.
$$
Since both $P_{\mathcal W}$ and $Q$ are orthogonal projections,
$$
 (P_{\cal W}-Q)^2
 =\begin{pmatrix}
 I-AA^*&0\\0&BB^*
 \end{pmatrix}
 =\begin{pmatrix}D&0\\0&BB^*\end{pmatrix}.
$$
Please noting  that $B^*B=I-A^*A$,
$$\sigma(BB^*)\setminus\{0\}=\sigma(B^*B)\setminus\{0\},$$ counting multiplicities. Similarly,
$$\sigma(I-AA^*)\cap (0,1)=\sigma(I-A^*A)\cap (0,1),$$ counting multiplicities.   Since $A$ is injective and $\dim\ker A^*=1$, the eigenvalue $1$ occurs once in $\sigma(D)$ and  $1\notin\sigma_p(B^*B)$.  Thus $1\in\sigma_p\big((P_{\cal W}-Q)^2\big)$ has multiplicity 1 and  $\lambda_j\in\sigma_{p}\big((P_{\cal W}-Q)^2\big)\cap(0,1)$ has multiplicities 2.
It is easy to verify
$$
 (P_{\cal W}-Q)(P_{\cal W}+Q-I)=-(P_{\cal W}+Q-I)(P_{\cal W}-Q), \text{ and }
 (P_{\cal W}+Q-I)^2=I-(P_{\cal W}-Q)^2.
$$
Let $E_{j}$ be the eigenspace  of $P-Q$ associated to the eigenvaluse $\lambda_j\in(0,1)$, the operator $(P+Q-I)|_{E_j}$ is invertible and commutes with the projections onto the eignspaces of $P-Q$ associated to the eigenvalues $+\sqrt{\lambda_j}$ and $-\sqrt{\lambda_j}$.  Hence these eigenvalues have equal multiplicities.  If $x\in\ker A^*$, then $x\in\mathcal W$, $Qx=0$, and $(P-Q)x=x$. Hence $1$ is the eigenvalue of $P_{\cal W}-Q$, with the multiplicity 1.  The spectral formula (\ref{eigenvalue}) follows, and therefore
$$
 \|P_{\cal W}-Q\|_{\mathcal S_p}^p
 =1+2\sum_j\lambda_j^{p/2}
 =2\left(1+\sum_j\lambda_j^{p/2}\right)-1
 =2\Tr(D^{p/2})-1.
$$
\end{proof}

Since $\|T_N\|=1$ and $0\le r<1$, the operator $I-rT_N$ is invertible. For simplicity, write
$$
 A_{N,r}=\varphi_r(T_N)=(rI-T_N)(I-rT_N)^{-1}, \text{ and } D_{N,r}=I-A_{N,r}A_{N,r}^*.
$$
 Moreover, since $T_N$ has no eigenvalues, we have $rI-T_N$ is injective and so is $A_{N,r}$.

Next, we claim that $\dim  \ker A_{N,r}^*=1$. In fact,
$$
 \ker A_{N,r}^*=\ker(rI-T_N^*).
$$
Writing $x=\sum\limits_{n=0}^\infty x_ne_n$, the equation $T_N^*x=rx$ is equivalent to
$$
 \alpha_{N,n}x_{n+1}=rx_n,\, n\ge0.
$$
After some direct calculations, $T_N^* K_r=rKr$, where
$$
 K_r=\sum_{n=0}^\infty
 \frac{r^n}{\alpha_{N,0}\alpha_{N,1}\cdots\alpha_{N,n-1}}e_n
 =\sum_{n=0}^\infty r^n\sqrt{\frac{h_{N,0}}{h_{N,n}}}\,e_n.
$$
By Stirling's formula, ${\color{red}h_{N,n}=\frac{(2N+n)!\,n!}{(N+n)!^2}}\to1$, and hence $K_r\in\ell^2$ for $r<1$.  Therefore
$$
 \dim\ker A_{N,r}^*=1.
$$

To continue, we need the following easy lemma.

\begin{lem}
\label{lem:defect-transform}
For $0\le r<1$,
\begin{eqnarray}\label{eq:Mobius}
 D_{N,r}
 =(1-r^2)(I-rT_N)^{-1}D_{N,0}(I-rT_N^*)^{-1}.
\end{eqnarray}
Consequently, for every $q>0$,
$$
 D_{N,r}\in\mathcal S_q\Leftrightarrow D_{N,0}\in\mathcal S_q \Leftrightarrow q>\frac12.
$$
\end{lem}

\begin{proof}It suffices to show (\ref{eq:Mobius}).
 Since $T_N$ commutes with $(I-rT_N)^{-1}$,
$$
 A_{N,r}=(I-rT_N)^{-1}(rI-T_N),\,  A_{N,r}^*=(rI-T_N^*)((I-rT_N)^{-1})^*.
$$
Therefore
\begin{align*}
 I-A_{N,r}A_{N,r}^*
 &=(I-rT_N)^{-1}\big[\big(I-rT_N)(I-rT_N^*)
 -(rI-T_N)(rI-T_N^*)]((I-rT_N)^{-1})^*\\
 &=(1-r^2)(I-rT_N)^{-1}(I-T_NT_N^*)((I-rT_N)^{-1})^*\\
 &=(1-r^2)(I-rT_N)^{-1}D_{N,0}(I-rT_N^*)^{-1}.
\end{align*}
Recall that
$
A_{N,0}=-T_N \text{ and }D_{N,0}=I-T_NT_N^*.
$
By Proposition \ref{prop:jacobi-model}, $D_{N,0} e_0=e_0$ and for $k\in\mathbb N$,
$$
D_{N,0}e_k=\frac{N^2}{(N+k)^2} e_k.
$$
It follows that $D_{N,0}\in \mathcal{S}_q$ if and only if $q> {1\over 2}$.
\end{proof}

Recall by (\ref{eq:6.2}),
$$
\mathbb U_r^* \mathfrak{C}_{\mathcal M_{N;\varphi_r}} \mathbb U_r=P-V_{r,w}PV_{r,w}^*,
$$
where $P$ is the projection onto $W_N$ and $A_{N,r}=PV_{r,w}|_{W_N}$. Since $A_{N,r}$ is injective and $\dim\ker A_{N,r}^*=1$, by Lemma \ref{lem:defect-transform}, $D_{N,r}$ is positive and compact. Applying Lemma \ref{lem:two-projections}, we have
$$
\|C_{\mathcal M_{N; \theta,\psi}}\|_{\mathcal S_p}^p=2\Tr(D_{N,r}^{p/2})-1.
$$
Moreover,
$$
C_{\mathcal M_{N; \theta,\psi}}\in\mathcal{S}_p\Leftrightarrow D_{N,r}\in\mathcal S_{p/2}\Leftrightarrow D_{N,0}\in \mathcal S_{p/2}\Leftrightarrow p>1.
$$ We summarize this argument  as the following result.
\begin{thm}
\label{thm:main}
Let $\theta,\psi$ be nonconstant inner functions, let $N\ge1$, and set
$$
 r=\rho_{\mathbb D}\big(\theta(0),\psi(0)\big).
$$
 Let $T_N$, $A_{N,r}$, and $D_{N,r}$ be defined as above.  Then
\begin{itemize}
\item[1)] the nonzero spectrum and nonzero singular values of $C_{\mathcal M_{N;\theta,\psi}}$ depend only on $N$ and $r$;
\item[2)] for every $0<p<\infty$,
$$
 C_{\mathcal M_{N;\theta,\psi}}\in\mathcal S_p
 \quad\Leftrightarrow\quad p>1;
$$
moreover, for $p>1$,
$$
 \|C_{\mathcal M_{N;\theta,\psi}}\|_{\mathcal S_p}^p
 =2\Tr(D_{N,r}^{p/2})-1.
$$
\end{itemize}
\end{thm}
\begin{cor}
Let $\theta$ and $\psi$ be two nonconstant inner functions, $\mathcal M_{N;\theta,\psi}=[(\theta-\psi)^N]$, then for $p>1$ we have the following estimation
$$
{1-r\over 1+r}\big(1+2N^p\zeta(p, N+1)\big)^{1\over p}\leq \|C_{\mathcal M_{N; \theta,\psi}}\|_{\mathcal S_p}\leq {1+r\over 1-r}\big(1+2N^p\zeta(p, N+1)\big)^{1\over p},
$$
where $r=\rho_{\mathbb D}\big(\theta(0),\psi(0)\big)$.
\end{cor}
\begin{proof} Combining with  Proposition \ref{Prop:trans} and
 Lemma \ref{lem:Composition},
the corollary comes from the classical norm formula; see e.g. \cite{Nordgren}, $$\|\mathcal{C}_{\varphi_r}\|=\|\mathcal C_{\varphi_r}^{-1}\|=\sqrt{1+r\over 1-r}.$$
\end{proof}

Finally, we will give the explicit formula of $\|C_{\mathcal M_{N; \theta,\psi}}\|_{\mathcal S_2}$. For convenience, setting
$$
 \varepsilon_{N,0}=1,
 \text{ and }
 \varepsilon_{N,n}=\frac{N}{N+n}\quad (n\in\mathbb N),
$$
one has that $D_{N,0}e_n=\varepsilon_{N,n}^2e_n$.  By the standard argument of Neumann series,
$$
 (I-rT_N)^{-1}e_\ell
 =\sum_{m=0}^\infty r^m
   \sqrt{\frac{h_{N,\ell+m}}{h_{N,\ell}}}\,e_{\ell+m}.
$$
Substituting this into Lemma \ref{lem:defect-transform}, we have, for $i,j\ge0$,
$$
 \langle D_{N,r}e_j,e_i\rangle
 =(1-r^2)\sqrt{h_{N,i}h_{N,j}}
 \sum_{\ell=0}^{\min\{i,j\}}
 \frac{\varepsilon_{N,\ell}^2r^{i+j-2\ell}}{h_{N,\ell}}.
$$
Since $D_{N,r}\ge0$, by the definition
\begin{align*}
 \Tr D_{N,r}
 =(1-r^2)\sum_{n=0}^\infty
   \varepsilon_{N,n}^2
   \left\|(I-rT_N)^{-1}e_n\right\|^2
 =(1-r^2)\sum_{n=0}^\infty\sum_{m=0}^\infty
   \varepsilon_{N,n}^2r^{2m}
   \frac{h_{N,n+m}}{h_{N,n}}.
\end{align*}
Therefore
$$
 \|C_{\mathcal M_{N;\theta,\psi}}\|_{\mathcal S_2}^2
 =2(1-r^2)\sum_{n=0}^\infty\sum_{m=0}^\infty
 \varepsilon_{N,n}^2r^{2m}\frac{h_{N,n+m}}{h_{N,n}}-1.
$$
At last, since
$$
 \sum_{m=0}^\infty r^{2m}\frac{h_{N,n+m}}{h_{N,n}}
 ={}_3F_2\!\left(
 \begin{matrix}1,n+1,n+2N+1\\n+N+1,n+N+1\end{matrix};r^2
 \right),
$$
where ${}_3F_2$ is the generalized hypergeometric function. This yields the explicit hypergeometric representation of the Hilbert-Schmidt norm:
$$
 \|C_{\mathcal M_{N;\theta,\psi}}\|_{\mathcal S_2}^2
 =2(1-r^2)\sum_{n=0}^\infty\varepsilon_{N,n}^2
 {}_3F_2\!\left(
 \begin{matrix}1,n+1,n+2N+1\\n+N+1,n+N+1\end{matrix};r^2
 \right)-1.
$$


\textbf{Acknowledgement.} The authors sincerely thank K. Guo and R. Yang for their encouragements and interests. The authors also  thank  R. Yang for valuable discussions on this topic. The second author is a member of LMNS, Fudan University.

\end{document}